\documentclass[hidelinks,onefignum,onetabnum]{siamart250211}

\usepackage{lipsum}
\usepackage{amsfonts}
\usepackage{graphicx}
\usepackage{epstopdf}
\ifpdf
  \DeclareGraphicsExtensions{.eps,.pdf,.png,.jpg}
\else
  \DeclareGraphicsExtensions{.eps}
\fi

\newsiamremark{remark}{Remark}
\newsiamthm{example}{Example}
\newsiamremark{hypothesis}{Hypothesis}
\crefname{hypothesis}{Hypothesis}{Hypotheses}
\newsiamthm{claim}{Claim}
\newsiamremark{fact}{Fact}
\crefname{fact}{Fact}{Facts}

\headers{Algebraic Varieties in Second Quantization}{Svala Sverrisdóttir}

\title{Algebraic Varieties in Second Quantization\thanks{Submitted to the editors May 23rd 2025.}}

\author{Svala Sverrisdóttir\thanks{Department of Mathematics, University of California at Berkeley
  (\email{svala@math.berkeley.edu})}}

\usepackage{amsopn}

\ifpdf
\hypersetup{
  pdftitle={Algebraic Varieties in Second Quantization},
  pdfauthor={Svala Sverrisdóttir}
}
\fi
\usepackage{graphicx}
\usepackage{caption}
\usepackage{mathtools}
\usepackage{enumerate}
\usepackage{verbatim}
\usepackage{tikz,tikz-cd,tikz-3dplot}
\usepackage{amssymb}
\usetikzlibrary{matrix}
\usetikzlibrary{arrows}
\usepackage{amsmath, amssymb, color}
\usepackage[noend]{algpseudocode}
\usepackage{caption}
\usepackage[normalem]{ulem}
\usepackage{subcaption}
\usepackage{multicol}
\usepackage{makecell}
\usepackage{array}
\usepackage{enumitem}
\usepackage{simpler-wick}

\newcommand{\PP}{\mathbb{P}}
\newcommand{\RR}{\mathbb{R}}

\newcommand{\CC}{\mathbb{C}}
\newcommand{\ZZ}{\mathbb{Z}}

\begin{document}
\maketitle

\begin{abstract} \noindent
We develop an algebraic geometric framework for Fock space coupled cluster theory in second quantization. In quantum chemistry, many-electron states are represented as elements of the exterior algebra. The fermionic creation and annihilation operators generate the Fermi–Dirac algebra, which can be realized as a Clifford algebra acting on the exterior algebra. We present a non-commutative Gröbner basis for the Fermi-Dirac algebra; offering an alternative proof of Wick's theorem, a fundamental result in
quantum field theory.
In coupled cluster theory, eigenpairs of the Schrödinger equation are approximated by a hierarchy of polynomial equations corresponding to different levels of truncation.
The coupled cluster exponential parameterization of quantum states gives rise to Fock space truncation varieties. This reveals well-known varieties, such as the Grassmannian, flag varieties and spinor varieties. We offer a detailed study of the truncation varieties, providing an explicit description of their defining equations and dimension. Furthermore, we classify all cases in which their coupled cluster degree coincides with the degree of the graph of the exponential parameterization -- most notably for singleton truncations such as CCD and for the Schubert like truncation varieties such as the Grassmannian, flag variety and the spinor variety. 
\end{abstract}

% REQUIRED

% REQUIRED

\section{Introduction}

Electronic structure theory is a powerful framework used in quantum chemistry to study the interactions of electrons in molecules \cite{helgaker2013molecular}. At the heart of it lies the time-independent electronic Schrödinger equation. In our setup we work over the fermionic Fock space $\mathcal{F}$, defined as the exterior algebra of a discretized Sobolov space spanned by $n$ one-particle basis functions. We are left with an innocent looking eigenvalue problem
$$
H\psi = \lambda \psi, \quad \psi \in \mathcal{F}
$$
where $H$ is a matrix representation of the \textit{Hamiltonian operator} on  $\mathcal{F}$. The size of the Hamiltonian matrix grows exponentially with $n$ and for large electronic systems, solving this Schrödinger equation is infeasible. Efficient and tractable numerical approximations, using e.g.\ perturbational and variational methods, are therefore essential for understanding the behavior of complex molecular structures. Within quantum chemistry, \textit{coupled cluster (CC) theory} is considered the gold standard for accurately approximating ground-state energies of weakly correlated systems near equilibrium. 

In earlier work, \cite{FSS}, with Fabian Faulstich and Bernd Sturmfels, we developed algebraic geometry for coupled cluster theory.
We introduced new algebraic objects, such as the truncation varieties -- projective varieties parameterizing feasible quantum states; and their CC degree, a complexity measure for solving the CC equations. The application of algebraic geometry to coupled cluster theory has been motivated by the works of scientists such as Piotr Piecuch, Karol Kowalski and Fabian Faulstich, who have long been interested in the existence of multiple solutions to the coupled cluster equations, e.g.\ \cite{kowalski2000complete}. For instance, \cite{piecuch2000search},  employs homotopy methods, identifying all CC solutions of a four-electron system in a standard quantum chemistry basis and \cite{faulstich2024coupled} provides an explicit upper bound for the CC degree in general. Earlier mathematical investigations of coupled cluster theory were carried out within a functional-analytic framework, beginning with Schneider in \cite{schneider2009analysis} and later extended by Rohwedder in \cite{rohwedder2013continuous}. More recently, Čsirik and Laestadius employed topological degree theory to analyze general coupled cluster variants in \cite{csirik2023coupled, csirik2023coupled2}.

The connection between quantum chemistry and algebraic geometry has turned out fruitful for both fields.
In algebraic geometry this has been explored by the author, Viktoriia Borovik and Bernd Sturmfels in \cite{borovik2025coupled} and by Evgeny Feigin in \cite{feigin2024birational, feigin2025birational}. The CC degree of the Grassmannian was identified with the degree of the graph of a given birational parametrization.
In contrast, on the quantum chemistry side, in \cite{sverrisdottir2024exploring}, the author and Fabian Faulstich explore all the CC solutions for LiH and H$_4$ in several symmetries. A recent paper by Fabian Faulstich, \cite{faulstich2024recent}, presents an in-depth overview of the latest mathematical developments in CC theory where the importance of finding all roots of the CC equations is explained. Algebraic geometry is central for understanding the root structure of the CC equations.

In first quantization one works over a discretized Hilbert space $\mathcal{H}$ equal to a  $d$-electron skew-symmetric product of a space spanned by $n$ one-particle basis functions.
In second quantization, see \cite{faulstich2024recent, faulstich2024coupled, szabo1989modern}, we extend from the $\binom{n}{d}$-dimensional space $\mathcal{H} \cong \wedge^d \RR^n$ to the Fock space, an exterior algebra $\mathcal{F} \cong \wedge \RR^n$ of dimension $2^n$. As an algebra, the Fock space offers richer structures, which are useful in applications, leading to faster computations using, for example, Wick's theorem. 
In this article, we extend the algebraic geometry of coupled cluster theory from first quantization, as was done in \cite{FSS},  to second quantization. 
% Our point of departure is the Fermi-Dirac algebra, a Clifford algebra generated by endomorphisms $a_p$ and $a_p^\dagger$ on $\mathcal{F}$ called the creation and annihilation operators. An element in the Fermi-Dirac algebra called the cluster operator is of special interest. Its $2^n$-dimensional matrix representation is used to define a parametrization of the quantum states, called the exponential map. We introduce a class of projective varieties in the space of binary tensors, $\PP^{2^n - 1}$, called the truncation varieties and are denoted $V_\sigma$. They are parametrized by restricting the domain of the exponential map. This procedure exposes as truncation varieties many widely studied structures such as flag varieties and spinor varieties. The CC equations are defined as a certain generalization of an eigenvalue problem over $V_\sigma$. The truncation varieties therefore represent the feasible quantum states in a given~truncation.
We now summarize the organization and contributions of this paper.

In Section \ref{se:FDalg}, we introduce the Fermi-Dirac algebra, a non-commutative algebra that can be realized as a Clifford algebra over $\mathcal{F}$. It is generated by the fermionic creation and annihilation operators -- exterior and interior products respectivly defined for each basis vector of $\RR^n$. 
The elements in the Fermi-Dirac algebra are endomorphisms on $\mathcal{F}$. 
In Theorem \ref{thm:wick}, we find a Gröbner basis for the two-sided ideal defining the Fermi-Dirac algebra using the diamond lemma for ring theory \cite[Theorem~1.2]{bergman1978diamond}. We explicitly describe the $2^n \times 2^n$ matrix representations of the elements in the Fermi-Dirac algebra using complete set matchings. This is an algebraic presentation of a foundational result in quantum field theory called Wick's theorem. 

In Section \ref{se:ExpParam}, we present the exponential parametrization -- a bijection with a polynomial inverse defined using an element in the Fermi-Dirac algebra, called the cluster operator. It expresses the quantum states in terms of cluster amplitudes. Theorem \ref{thm:psiparam} gives an explicit formula, in terms of even set partitions, for all the coordinates of both the forward and inverse map. The extension from first to second quantization enables us to fully describe the exponential parametrization of quantum states using matrix representations of endomorphisms of the Fock space $\mathcal{F}$. 

Section \ref{se:TruncVars} introduces the Fock space truncation varieties $V_\sigma$ -- a family of projective varieties in the space of binary tensors, $\PP^{2^n - 1}$. They are parametrized by restricting the domain of the exponential map and therefore represent the feasible quantum states in a given~truncation.
The restriction of $V_\sigma$ to an affine chart is a complete intersection, defined by some coordinates of the
inverse exponential map. The truncation varieties introduced in \cite[Section~3]{FSS}, embedded into $\PP^{2^n - 1}$, are an example of Fock space truncation varieties, revealing, for example, the Grassmannian.
Proposition \ref{prop:parthole} features an isomorphism between truncation varieties corresponding to $(d,n)$ and $(n - d,n)$ using the particle-hole symmetry. Theorem \ref{thm:linear} identifies all linear Fock space truncation varieties. Finally, Theorem \ref{thm:eomcc} describes the relation between our Fock space truncation varieties and the ones in EOM-CC.

Section \ref{se:FlagsSpinors} centers around two well-known varieties in the space of binary tensors, specifically, the flag variety and the spinor variety. Theorem \ref{thm:flag} reveals the partial flag variety $\mathcal{F}l(d - 1, d, d + 1, n)$ as a truncation variety. In particular, the Flag variety parallels the Grassmannian as the truncation variety for CCS in the case of ionization and electron attachment. The exponential map displays a unique parametrization of the flag variety. In Theorem \ref{thm:spin}, we show that the spinor variety, parametrized by Pfaffians of a skew-symmetric matrix, also appears as a Fock space truncation variety. 

Section \ref{se:CCeqs} studies the CC equations in second quantization. They are defined as a certain generalization of an eigenvalue problem over a truncation variety $V_\sigma$. The number of solutions is finite for a generic Hamiltonian $H$. This number is an invariant of the truncation variety, called the CC degree. It serves as a complexity measure for solving the CC equations.
In \cite{borovik2025coupled}, we identified the CC degree of the Grassmannian with the degree of the graph of its parameterization. In Theorem \ref{thm:degofgraph} we classify all truncation varieties with this property. Two such varieties are the flag variety and the spinor variety. It would be of interest to investigate toric degenerations that lift to the graph of the parametrizations introduced in Section~\ref{se:FlagsSpinors}. This could reveal a description or an explicit formula for their CC degrees. 

Section~\ref{sec:numsols} explores the computational aspects of solving the coupled cluster equations using algebraic methods, such as, parameter homotopy and monodromy.
We numerically calculate the CC degrees of $\mathcal{F}l(d - 1, d, d + 1, n)$ and the spinor variety in small cases using \texttt{HomotopyContinuation.jl} \cite{breiding2018homotopycontinuation}. We explore the complexity of computing the solution spectrum of the CCSD equations when working with ionized electronic systems. We conclude in Example 7.3 with the case of 3-electron systems.

\section{The Fermi-Dirac algebra}\label{se:FDalg}

We fix a positive integer $n$. In second quantization, we work in the \textit{Fock space}, denoted by $\mathcal{F}$. It is the exterior algebra of a discretized Hilbert space of sufficiently smooth real functions on the product of $\RR^3$ with the spin set $\{\pm \frac{1}2{}\}$. This Hilbert space is called the spin-orbital space and it is the Sobolev space $H^1(\{\pm \frac{1}2{}\} \times \RR^3 )$. The discretization is an $n$-dimensional approximation of the spin-orbital space, derived by choosing a suitable finite basis of functions. Our $n$ basis functions are called the \textit{molecular spin orbitals}. They are linear combinations of \textit{atomic spin orbitals} -- harmonic-like functions describing the charge distribution around an atom's nucleus. See \cite{helgaker2013molecular} for a broad introduction into electronic structure theory and \cite[Section~3]{faulstich2024coupled} for a detailed description of second quantization. In this paper, we work within the spin-generalized setting and do not impose any spin restrictions.
We note that the Fock space is isomorphic to the exterior algebra of $\mathbb{R}^n$, that is, $\mathcal{F} \cong \wedge \mathbb{R}^n$. Throughout this paper, we will work within this exterior algebra, as well as in its complexification and projectivization, $\wedge \mathbb{C}^n$ and $\mathbb{P}(\wedge \mathbb{C}^n) \cong \mathbb{P}^{2^n - 1}$.

The standard basis vector $e_p$ of $\RR^n$ corresponds to the $p$-th molecular spin orbital. The basis vectors of $\mathcal{F}$ are denoted by $e_J = e_{j_1} \wedge \cdots \wedge e_{j_k}$ with index sets $J = \{j_1, \dots, j_k\} \subseteq [n] = \{1,2,\dots,n\}$. The size of the set $k = |J|$ can vary from $0$ to $n$. The basis is ordered in a reverse lexicographic order on the index sets $J$. That is:
$$
\varnothing \le 1 \le 2 \le 12 \le 3 \le 13 \le 23 \le 123 \le 4 \le 14 \le 24 \le 124 \le 34 \le 134 \le 234 \le \cdots.
$$
In quantum chemistry, the basis vectors of $\mathcal{F}$ are referred to as \textit{occupation number states}, a term we will also use here. They correspond to a configuration of particles on the molecular orbitals. Pauli's exclusion principle states that two particles cannot simultaneously occupy the same spin-orbital. Hence, the basis vectors describe every configuration of particles on our electronic system. We say the $p$-th spin orbital is \textit{occupied} in an occupation number state $e_J$, if $p \in J$; otherwise it is \textit{unoccupied}. 

The quantum states are a linear combination of the occupation number states
$$
\psi = \sum_{J \subseteq [n]} \psi_J e_J.
$$
We call the coordinates, $\psi_J$, \textit{Plücker coordinates}, highlighting the connection to the Grassmannian and the flag variety. The coefficient $\psi_J$ corresponds to a  minor of an $n \times n$ matrix $\Theta$, taking rows $1, \dots, k$ and columns $j_1, \dots j_k$, where $k = |J|$. These are the Plücker coordinates of the matrix $\Theta$, parameterizing the complete flag spanned by the rows of $\Theta$, \cite[Chapter~14.1]{miller2005combinatorial}.

We fix another positive integer $d \le n$, enumerating the base number of electrons in an electronic system. We call $e_{[d]}$ the \textit{reference state}. The quantum states can be written in terms of the {\em configuration interaction coefficients} $c_{I,B}$, where $I \subseteq [ d ]$ and $B \subseteq [ n ]\backslash [ d ]$. The correspondence between the two coordinate systems is defined as:
\begin{equation}\label{eq:coordcorrespond}
    \psi_J = c_{I, B} \quad\text{where}\quad J = ([d] \backslash I) \cup B \quad\text{and}\quad I = [d] \backslash J, \,\, B = J \backslash [d].
\end{equation}
See \cite[Equation~1]{FSS} for the equivalent definition in first quantization. The configuration interaction coefficient of the reference state is $c_{0,0}$. The number of configuration interaction coefficients is 
$$
\sum_{m,\ell = 0}^{d, n - d} \binom{d}{m}\binom{n - d}{\ell} = \sum_{k = m + \ell = 0}^n \binom{n}{k} = 2^n
$$
so the two coordinate sets have the same cardinality.

In the exterior algebra, $\mathcal{F} \cong \wedge \RR^n$, we define the \textit{exterior and interior product} for each basis element $e_p \in \RR^n$ as the following endomorphisms:
\begin{equation}
\label{eq:CreandAni}
    \begin{aligned}
        &a_p^\dag: \mathcal{F} \to \mathcal{F}, \quad \psi \,\mapsto\, e_p \wedge \psi \,=\, \sum_{I \subseteq [n]} \psi_I  e_p \wedge e_I\\
        &a_p: \mathcal{F} \to \mathcal{F}, \quad \psi \,\mapsto\, e_p~\lrcorner~\psi \,=\, \sum_{I \subseteq [n]} \psi_I  e_p~\lrcorner~e_I.
    \end{aligned}
\end{equation}
Here $\lrcorner$ denotes the dual operation of the wedge product $\wedge$, called the hook product \cite[Section~3.6]{gallier2020differential}. One can define the hook product using the particle hole symmetry in the following way: 
$$
e_p~\lrcorner~\psi = (e_p \wedge \psi^\dagger)^\dagger.
$$
The operators in (\ref{eq:CreandAni}) are called the \textit{creation} and \textit{annihilation operators}, respectively. As the name suggests, the physical meaning of them is to create or annihilate a particle in the $p$-th molecular spin orbital.

\begin{remark}[Jordan-Wigner transformation]\label{re:matrixrep}
We describe a $2^n$-dimensional matrix representation of endomorphisms on the Fock space $\mathcal{F}$, using the Jordan-Wigner transformation \cite{jordan1928pauli}. For an English reference on the transformation we point the reader to \cite{nielsen2005fermionic}. An unoccupied orbital is represented by 
$\begin{bmatrix}
    1 & 0
\end{bmatrix}^T$, and an occupied orbital by $\begin{bmatrix}
    0 & 1
\end{bmatrix}^T$. An occupational number state $e_J$ can be represented by a tensor product of the $2$-dimensional occupied and unoccupied vectors. For example, the reference state is represented by
$$
e_{[d]} \equiv \underbrace{\begin{bmatrix}
    0 \\ 1
\end{bmatrix} \otimes \cdots \otimes \begin{bmatrix}
    0 \\ 1
\end{bmatrix} }_{d \text{ times}} \otimes \underbrace{\begin{bmatrix}
    1 \\ 0
\end{bmatrix} \otimes \cdots \otimes \begin{bmatrix}
    1 \\ 0
\end{bmatrix}}_{n - d \text{ times}}.
$$
The creation and annihilation operators can be represented by $2^n \times 2^n$ matrices as follows. Let
$$
\sigma_z = \begin{bmatrix}
    1 & 0\\
    0 & -1
\end{bmatrix}, \quad a = \begin{bmatrix}
    0 & 1\\
    0 & 0
\end{bmatrix}
$$
be the $2 \times 2$ \textit{Pauli-z matrix} and \textit{annihilation matrix}.
The annihilation matrix, $a$, turns an occupied orbit into an unoccupied orbit by left matrix multiplication. Its adjoint, $a^\dagger$, is the \textit{creation matrix}, turning unoccupied orbits into occupied orbits by left multiplication. The $2^n \times 2^n$ creation and annihilation matrices for the $p$-th orbital are defined as tensor products of $2 \times 2$ matrices:
\begin{equation}
\label{eq:CreationAndAnnihilation}
a_p^\dagger
=
\underbrace{\sigma_z \otimes ... \otimes \sigma_z}_{p-1~{\rm times}} 
\otimes\;
a^\dagger
\otimes 
\underbrace{I_2 \otimes ... \otimes I_2}_{n-p~{\rm times}},
\quad
a_p
=
\underbrace{\sigma_z \otimes ... \otimes \sigma_z}_{p-1~{\rm times}} 
\otimes\;
a
\otimes 
\underbrace{I_2 \otimes ... \otimes I_2}_{n-p~{\rm times}}.
\end{equation}
    
\end{remark}

A \textit{Clifford algebra} on the free associative algebra $\CC \langle x_1, \dots x_k \rangle$  is a quotient defined by a two-sided ideal of the form
$$
\langle \, x_i x_j + x_j x_i - q(x_i, x_j) 1 \, : \, 1 \le i,j \le k \, \rangle
$$
where $q$ is a quadratic form on $\CC \langle x_1, \dots x_k \rangle$ and $1$ is the multiplicative identity. For an in-depth introduction to Clifford algebras and their connection to physics, we refer to Chevalley's book \cite{chevalley1995algebraic}. They arise in our quantum chemistry formulation~as follows:

\begin{proposition}\label{prop:cliff}
The creation and annihilation operators span the Fermi-Dirac algebra, which is a special Clifford algebra on $\CC\langle a_1, \dots, a_n, a_1^\dagger, \dots a_n^\dagger \rangle$, associated with the quadratic form
\begin{equation*}
    Q = \begin{bmatrix}
    0_{n} & I_{n}\\
    I_{n} & 0_{n}\\
\end{bmatrix}.
\end{equation*}
\end{proposition}

\begin{proof}
The creation/annihilation operators fulfill the anti-commutator relations
\begin{equation}\label{eq:fermirels}
    a_p^\dag a_q^\dag + a_q^\dag a_p^\dag = 0, \quad a_pa_q + a_qa_p = 0, \quad a_p^\dag a_q + a_q a_p^\dag = \delta_{p,q}.
\end{equation}
This means they span the Fermi-Dirac algebra.
\end{proof}

\begin{corollary}\label{cor:crannilpotent}
    The creation and annihilation operators, respectively, span an isotropic vector space of dimension $n$, with respect to quadratic form $Q$. Hence, they are nilpotent of order~$1$.
\end{corollary}

\begin{proof}
    Since $Q$ has zero blocks on the diagonal, we see that the vector spaces spanned by the creation and annihilation operators are isotropic, respectively.
\end{proof}

\begin{remark}[Exterior algebras in the Fermi-Dirac algebra]\label{re:exteriorsofFD}
    We look at a subalgebra $\mathcal{S}$, of the Fermi-Dirac algebra, generated by the $n$ creation and annihilation operators, $a_1,\dots,a_d, a_{d + 1}^\dagger, \dots a_n^\dagger$. This is a Clifford algebra associated to a quadratic form, defined by a submatrix of $Q$, equal to the zero matrix. Hence, $\mathcal{S}$ is isomorphic to the exterior power $\wedge \CC^n$. A word $\Omega \in \mathcal{S}$ of even length, $\ell(\Omega)$, commutes with any element in $\mathcal{S}$.  
    More specifically, words $\Omega$ and $\Omega'$
    in $\mathcal{S}$ fulfill the following commutator and anti-commutator relations
    $$
    [\,\Omega , \Omega'\,] =0 ,\text{ if } \ell(\Omega) \text{ or } \ell(\Omega') \text{ is even}, \quad [\,\Omega, \Omega'\,]_+ =0 ,\text{ if } \ell(\Omega) \text{ and } \ell(\Omega') \text{ are odd}.
    $$
    More specifically, even strings in $\mathcal{S}$ commute with elements in $\mathcal{S}$.
    
    There exist $2^n$ subalgebras of the Fermi-Dirac algebra isomorphic to $\wedge \CC^n$. We construct such a subalgebra by choosing either a creation or annihilation operator
    for each index $i \in [n]$.
\end{remark}

We fix the following variable order on the creation and annihilation operators:
    $$
    a_n > a_{n - 1} > \cdots > a_1 > a_1^\dagger > a_2^\dagger > \cdots > a_n^\dagger.
    $$
The degree lexicographic order on monomials in the free associative algebra is defined as follows: monomial $w = w_1 \cdots w_\ell$ of degree $\ell$ is less than monomial $w' = w'_1 \cdots w'_p$ of degree $p$ if $\ell < p$ or $\ell = p$ and for some $k$ both $w_k < w_k'$ and  $w_{i} = w'_{i}$ for all $i < k$.
For example
$$
a_1^\dagger < a_2 < a_1^\dagger a_2^\dagger < a_1a_2^\dagger < a_1^{\dagger 2} < a_1a_2 < a_2^2 < a_1^3 < a_1a_2^2.
$$
For a detailed review on non-commutative algebras and their Gröbner bases, we refer to Bergman's foundational paper on the diamond lemma for ring theory, \cite{bergman1978diamond}, and a paper by Teo Mora on commutative and non-commutative Gröbner bases, \cite{mora1994introduction}.

\begin{theorem}[Wick's theorem \cite{wick1950evaluation}]\label{thm:wick}
The relations in (\ref{eq:fermirels}), generating the two-sided ideal, defining the Fermi-Dirac algebra, form a Gröbner basis with respect to degree lexicographic order. 
\end{theorem}

\begin{proof}
    First, we let $G$ be the set of relations in (\ref{eq:fermirels}), that is
    $$
    G = \{ \underline{a_i a_j}  + a_ja_i \, : \, i \ge j \}\, \cup \,\{ \underline{a_i^\dagger a_j^\dagger}  + a_j^\dagger a_i^\dagger \, : \, i \le j \} \,\cup\, \{ \underline{a_i a_j^\dagger}  + a_j^\dagger a_i - \delta_{ij} \, : \, i,j \}.
    $$
    The leading terms are underlined. We define a \textit{critical pair} of $G$ to be a pair of leading terms $(l(f), l(f'))$ of elements $f,f'$ in $G$ such that $l(f) u = v l(f')$ for some variables $u$ and $v$ in $\CC\langle a_1, \dots, a_n, a_1^\dagger, \dots a_n^\dagger \rangle$. The \textit{S-elements} are the polynomials $fu - vf'$ where $l(f)$ and $l(f')$ form a critical pair \cite[Section~5.3]{mora1994introduction}. By the diamond lemma for ring theory \cite[Theorem~1.2]{bergman1978diamond} we only need to check that the S-elements of $G$ have a weak Gröbner representation, see \cite[Section~5.3]{mora1994introduction}. 
    First we find the critical pairs of $G$. They are:
    \begin{align*}
        &(a_i a_j, a_ja_k) \text{ where } i \ge j  \ge k, \,\,  (a_i^\dagger a_j^\dagger, a_j^\dagger a_k^\dagger) \text{ where }i \le j \le k,\\ &(a_i a_j, a_ja_k^\dagger) \text{ where } i \ge j, \,\,  (a_i a_j^\dagger, a_j^\dagger a_k^\dagger)\text{ where } j \le k.
    \end{align*}
    The S-elements are
    \begin{align*}
        &(a_ia_j + a_ja_i)a_k - a_i(a_ja_k + a_ka_j) = a_ja_ia_k - a_ia_ka_j\\
        &(a_i^\dagger a_j^\dagger + a_j^\dagger a_i^\dagger )a_k^\dagger - a_i^\dagger(a_j^\dagger a_k^\dagger + a_k^\dagger a_j^\dagger ) = a_j^\dagger a_i^\dagger a_k^\dagger - a_i^\dagger a_k^\dagger a_j^\dagger \\
        &(a_ia_j + a_ja_i)a_k^\dagger - a_i(a_ja_k^\dagger + a_k^\dagger a_j - \delta_{ik}) = a_ja_ia_k^\dagger - a_ia_k^\dagger a_j + \delta_{ik} a_i\\
        &(a_ia_j^\dagger + a_j^\dagger a_i - \delta_{ij})a_k^\dagger - a_i(a_j^\dagger a_k^\dagger + a_k^\dagger a_j^\dagger) = a_j^\dagger a_ia_k^\dagger - a_ia_k^\dagger a_j^\dagger - \delta_{ij} a_k^\dagger.
    \end{align*}
    One can check that all these polynomials have a weak Gröbner representation.
\end{proof}

The standard monomials of $\langle G \rangle$ are the following monomials
\begin{equation}\label{eq:standardMons}
    a_B^\dagger a_I = a_{b_\ell}^\dagger \cdots a_{b_1}^\dagger a_{i_1} \cdots a_{i_m}.
\end{equation}
In quantum chemistry these monomials are commonly referred to as normal ordered strings.
The Fermi-Dirac algebra is isomorphic to a vector space over $\CC$, generated by the standard monomials. There is a standard monomial for each pair of subsets $B \subseteq [n]$ and $I \subseteq [n]$, so the dimension of the Fermi-Dirac algebra as a vector space is $2^{2n}$.
Each element in the Fermi-Dirac algebra can be written in terms of the standard monomials. This is called the \textit{standard representation} of an element.

    We want to find the standard representation of the word 
    $\Omega = \omega_1 \cdots \omega_k$ in the Fermi-Dirac algebra. First, we set $I$ as the index set of the annihilation operators and $B$ as the index set of the creation operators in $\Omega$. Without loss of generality, we assume that $m= |I| \le |B| = \ell$. 
    We define the ordered term of $\Omega$~as 
    $$
    o(\Omega) = \operatorname{sign}(\Omega) a_{B}^\dagger a_{I}
    $$
    where $\operatorname{sign}(\Omega)$ is the sign of the permutation $\Omega \mapsto a_{B}^\dagger a_{I}$.
    This is the standard representation of $\Omega$ in the exterior algebra, that is, if all the creation and annihilation operators would anti-commute. However, we have relations 
    $
    a_p^\dag a_q + a_q a_p^\dag = \delta_{p,q}    
    $
    so the operators do not all anti-commute. We define the contraction of two operators as
    $$
    \wick{\c1 \omega \c1 \omega' } = \omega \omega' - o(\omega \omega ').
    $$
    The contraction is zero unless $\omega = a_p$ and $\omega' = a_p^{\dagger}$ for some $p$, then it is $1$.
    The following proposition describes the standard representation of words in the Fermi-Dirac algebra. Specifically, the constant term is of importance, as it describes elements in the matrix representation.
    % The standard monomials are endomorphisms on $\mathcal{F}$  send the $1 \equiv e_\varnothing \in \mathcal{F}$ 
    \begin{proposition}\label{prop:standrep}
        The standard representation of a word $\Omega = \omega_1 \cdots \omega_k$ in the Fermi-Dirac algebra is 
    \begin{align*}
        &\Omega \,\, = \,\, o(\Omega) + \sum_{\substack {i \in I, \,  b \in B\\ i < b}} (-1)^{i + b - 1}o(\omega_1 \cdots \wick{\c1 \omega_i \cdots \c1 \omega_b} \cdots \omega_k)\\
        &+ \sum_{\substack {i < j \in I, \,  b \ne c \in B\\ i < b, j < c}}  (-1)^{i + j + b + c - 2}(-1)^{\#\text{cross}} o(\omega_1 \cdots \wick{\c1 \omega_i \cdots \c2 \omega_j \cdots \c1 \omega_b \cdots \c2 \omega_c} \cdots \omega_k) + \,\, \cdots \,\, \\
        &+\sum_{\substack{i_1 < \cdots < i_m \in I, \\b_1 \ne \cdots \ne b_m \in B,\\ i_1 < b_1, \dots, i_m < b_m}} \hspace{-1.5em} (-1)^{i_1 + b_1 + \cdots + i_m + b_m - m}(-1)^{\#\text{cross}} o(\wick{\c2 \omega_{i_1} \c3 \omega_{i_2} \cdots \c1 \omega_{i_j} \cdots \omega_c \c1 \omega_{b_j} \cdots \c2 \omega_{b_1} \cdots \c3 \omega_{b_2}})
    \end{align*}
    where $\#$cross is the number of crossings between contractions. Each summand corresponds to a specific number of contractions.
    \end{proposition}

A \textit{complete matching} of the set $[2m]$ is defined to be a set of $m$ pairs $(i,j)$ where $i < j$ and each element in $[2m]$ appears exactly once. We say the element $i$ is the \textit{lefthand endpoint} of pair $(i,j)$ and the element $j$ is the \textit{righthand endpoint}, see \cite[Section~1]{chen2007crossings}. We notice that the standard representation of $\Omega$ sums over complete matchings of subsets of $[2m]$ where the lefthand endpoints are in $I$ and the righthand endpoints are in $B$.

\begin{example}[$a_ia_p^\dagger a_j a_q^\dagger$] We write the word $a_ia_p^\dagger a_j a_q^\dagger$ in standard representation:
\begin{align*}
    a_ia_p^\dagger a_j a_q^\dagger &= -a_p^\dagger a_q^\dagger a_ia_j + o(\wick{ \c1 a_i \c1 a_p^\dagger} a_j a_q^\dagger) + o(\wick{ \c1 a_i a_p^\dagger a_j \c1 a_q^\dagger}) +  o(\wick{ a_i a_p^\dagger \c1 a_j \c1 a_q^\dagger}) + \wick{ \c1 a_i \c1 a_p^\dagger \c1 a_j \c1 a_q^\dagger}\\
    &= -a_p^\dagger a_q^\dagger a_ia_j - \delta_{ip}a_q^\dagger a_j + \delta_{iq}a_p^\dagger a_j - \delta_{jq}a_p^\dagger a_i + \delta_{ip}\delta_{jq}.
\end{align*}    
\end{example}

\begin{remark}\label{re:constSR}
    The standard representation of a word $\Omega$ in the Fermi-Dirac algebra has a constant term only if there is an equal number of annihilation operators and creation operators, that is $|I| = |B| = m$. In that case, the parity of the first sign of the constant term in Proposition \ref{prop:standrep} is
    $$
    {i_1 + b_1 \cdots + i_m + b_m - m} = \sum_{i = 1}^{2m} i - m = \frac{(2m + 1)2m}{2} - m = 2m^2.
    $$
    Hence it is even, and we get a positive sign.
    The constant term of the standard representation of $\Omega$ is then equal to   
    $$
    \sum (-1)^{\#\text{cross}} \wick{\c2 \omega_{i_1} \c3 \omega_{i_2} \cdots \c1 \omega_{i_j} \cdots \c1 \omega_{b_j} \cdots \c2 \omega_{b_1} \cdots \c3 \omega_{b_2}}
    $$
    where we sum over all complete matchings of $[2m]$ where the lefthand endpoints are in $I$ and the righthand endpoints are in $B$. We note that a fully contracted word $\Omega$, is a product of $\delta$-functions. Therefore only complete matchings of the form $\{\,(a_p, a_p^\dagger) \, : \, i \,\}$ produce a non-zero contraction. There only exists such a matching if the annihilation indices $I$ are equal to the creation indices $B$.
\end{remark}

    A word $\Omega$ in the Fermi-Dirac algebra has a $2^n \times 2^n$ matrix representation as an endomorphism of $\mathcal{F}$. We can find it by multiplying together the matrices of the creation and annihilation operators in $\Omega$. This procedure involves many matrix multiplications for matrices of size $2^n \times 2^n$. For a large number $n$, this computation is not feasible. 

\begin{remark}\label{re:marepFD}
Alternatively we can define the matrix representation of $\Omega$ by evaluating at occupation number states. That is
$
(e_I^\dagger \Omega e_J)_{I,J \subseteq [n]}.
$
We notice that the occupation number state $e_I$ can be written as
$e_I = a_{i_1}^\dagger\cdots a_{i_m}^\dagger 1$
so 
$$
e_I^\dagger \Omega e_J = 1^\dagger a_{i_m} \cdots a_{i_1} \Omega a_{j_1}^\dagger \cdots a_{j_\ell}^\dagger 1.
$$
For every standard monomial $a_{B}^\dagger a_{I}$ we see that
$
1^\dagger a_{B}^\dagger a_{I} 1 = 0 
$. Using the standard representation of $a_{i_m} \cdots a_{i_1} \Omega a_{j_1}^\dagger \cdots a_{j_\ell}^\dagger$ we obtain the following description of the elements in the matrix representation
$$
e_I^\dagger \Omega e_J  = \text{constant term in standard rep. of } a_{i_m} \cdots a_{i_1} \Omega a_{j_1}^\dagger \cdots a_{j_\ell}^\dagger.
$$
By Remark \ref{re:constSR} the constant term can be calculated by finding the number of crossings in complete matchings $\{(a_p, a_p^\dagger)  :  p  \}$, if they exists; otherwise the constant term is zero. 
\end{remark}

Electronic operators such as the Hamiltonian and the cluster operator -- introduced in the following section -- are elements of the Fermi–Dirac algebra. By employing Remark~\ref{re:marepFD}, the coupled cluster equations, seen in Section~\ref{se:CCeqs}, can be derived through straightforward combinatorial manipulations. In contrast, the first-quantized operators are represented by $\binom{n}{d} \times \binom{n}{d}$ matrices. Then, deriving approximation schemes like the coupled cluster equations requires prohibitively large matrix multiplications. This highlights the advantage of the second-quantized formalism, where corollaries of Wick’s theorem enable efficient derivations. This framework for coupled cluster theory was first established in the foundational 1966 paper~\cite{vcivzek1966correlation}.

\section{Exponential parameterization}\label{se:ExpParam}

This section generalizes \cite[Section~2]{FSS} from first quantization to second quantization. There we defined a parametrization of the quantum states, called the exponential parametrization, writing them as a column of the exponential matrix of a $\binom{n}{d} \times \binom{n}{d}$ matrix $T(t)$. We  assumed the particle number was fixed, denoted fixed-N, and worked within the $d$-th grading of the Fock space, denoted $\mathcal{H} \cong \wedge^d \RR^n$.
By extending to the whole Fock space we are able to describe the exponential parametrization explicitly using an element of the Fermi-Dirac algebra. Now we do not have to make the assumption that the number of particles is fixed. Indeed, this framework also captures cases of ionization and electron attachment. Our goal is to build a general algebraic framework for coupled cluster theory encapsulating fixed-N CC as well as Fock space coupled cluster (FSCC) \cite{eliav2021relativistic,kaldor1991fock, lindgren1987connectivity} and equations of motion coupled cluster (EOM-CC) \cite{shavitt2009many}, see  for example Theorem \ref{thm:eomcc}.  

We take another copy of the exterior algebra $\mathcal{V} = \wedge \RR^n$ whose elements $t$ are called \textit{cluster amplitudes}. Their coordinates $t_{I,B}$ are indexed by $I \subseteq [ d ]$ and $B \subseteq [ n ] \backslash [ d ]$. We also introduce alternative coordinates $x_J$ obtained from $t_{I,B}$ by the correspondence $J = ([d] \backslash I) \cup B$. 
We define the \textit{annihilation and creation level} of the coordinates for the exterior algebras as:
$$
a(\psi_J) = a(x_J) = a(J) = |[d] \backslash J |\,\,\,\text{and}\,\,\, c(\psi_J) = a(x_J) = c(J) = |J \backslash [d]|.
$$
The \textit{level} of coordinates $\psi_J$, $x_J$ is defined to be the sum of the annihilation and creation level. In terms of coordinates $t_{I,B}$ and $c_{I,B}$, the annihilation level is equal to $|I|$, the creation level is $|B|$ and the level of the coordinate is $|I| + |B|$.

We recall the Fermi-Dirac algebra has a non-commutative Gröbner basis defined by the anti-commutator relations (\ref{eq:fermirels}). Hence every element in the algebra can be written in terms of the standard monomials (\ref{eq:standardMons}).
We define the \textit{cluster operator} to be a sum, with unknown coefficients, of the standard monomials with positive degree whose annihilation operators are indexed by subsets of $[d]$ and creation operators are indexed by subsets of $[n] \backslash [d]$. In symbols: 
\begin{equation}
T(t) = \sum_{\substack{ I \subseteq [d], \,  B \subseteq [n] \backslash [d] \\ |I| +  |B| > 0}} t_{I,B} a_{b_\ell}^\dagger...a_{b_1}^\dagger
a_{i_1}...a_{i_m} =  \sum_{\substack{ I \subseteq [d], \,  B \subseteq [n] \backslash [d] \\ |I| +  |B| > 0}} t_{I,B} a_B^\dagger
a_I.
\end{equation}
The cluster operator is an element of the Fermi-Dirac algebra. In fact, $T(t)$ is a general element of the subalgebra $\mathcal{S} \cong \wedge \CC^n$ of the Fermi-Dirac algebra, defined in Remark \ref{re:exteriorsofFD}.

The matrix representation of $T(t)$ is called the \textit{cluster matrix}. 
By Corollary \ref{cor:crannilpotent} and the pigeonhole principle, the cluster matrix $T(t)$ is nilpotent of order $n$, hence $T(t)^{n + 1} = 0$. We can therefore define the matrix exponential of $T(t)$ as a finite sum
$$
\exp(T(t)) = \sum_{k = 0}^n \frac{1}{k!} T(t)^k.
$$
See \cite[Section~3]{faulstich2024coupled} for a more detailed overview of $T(t)$ and its properties. We define the \textit{exponential parametrization} to be the  map
\begin{equation}\label{eq:expparam}
    \mathcal{V} \rightarrow \mathcal{F}, \,\,\,  t \mapsto \psi =  \exp(T(t))e_{[d]}.
\end{equation}
This map defines a parametrization of the quantum states in the Fock space $\mathcal{F}$.  This parametrization is used when employing the Fock space coupled cluster method -- a widely studied theory for handling ionized electronic systems.

\begin{example}[$d = 2$, $n = 4$]
We look at the case when we have four atomic orbitals and two base particles. The cluster matrix is the $16 \times 16$ matrix of the form:

\begin{tiny}
$$
   T(t) =
\left [ \arraycolsep=1.9pt
\begin{array}{cccccccccccccccc}
   0 & t_{1,0} & t_{2,0} & t_{12,0} & 0 & 0 & 0 & 0 & 0 & 0 & 0 & 0 & 0 & 0 & 0 & 0 \\
0 & 0 & 0 & -t_{2,0} & 0 & 0 & 0 & 0 & 0 & 0 & 0 & 0 & 0 & 0 & 0 & 0 \\
0 & 0 & 0 & t_{1, 0} & 0 & 0 & 0 & 0 & 0 & 0 & 0 & 0 & 0 & 0 & 0 & 0 \\
0 & 0 & 0 & 0 & 0 & 0 & 0 & 0 & 0 & 0 & 0 & 0 & 0 & 0 & 0 & 0 \\
t_{0,3} & t_{1,3} & t_{2,3} & t_{12,3} & 0 & t_{1,0} & t_{2,0} & t_{12,0} & 0 & 0 & 0 & 0 & 0 & 0 & 0 & 0 \\
0 & -t_{0,3} & 0 & t_{2 ,3} & 0 & 0 & 0 & -t_{2,0} & 0 & 0 & 0 & 0 & 0 & 0 & 0 & 0 \\
0 & 0 & -t_{0,3} & -t_{1 ,3} & 0 & 0 & 0 & t_{1,0} & 0 & 0 & 0 & 0 & 0 & 0 & 0 & 0 \\
0 & 0 & 0 & t_{0,3} & 0 & 0 & 0 & 0 & 0 & 0 & 0 & 0 & 0 & 0 & 0 & 0 \\
t_{0,4} & t_{1 ,4} & t_{2,4} & t_{12 ,4} & 0 & 0 & 0 & 0 & 0 & t_{1,0} & t_{2,0} & t_{12,0} & 0 & 0 & 0 & 0 \\
0 & -t_{0,4} & 0 & t_{2,4} & 0 & 0 & 0 & 0 & 0 & 0 & 0 & -t_{2,0} & 0 & 0 & 0 & 0 \\
0 & 0 & -t_{0,4} & -t_{1 ,4}  & 0 & 0 & 0 & 0 & 0 & 0 & 0 & t_{1,0} & 0 & 0 & 0 & 0 \\
0 & 0 & 0 & t_{0,4} & 0 & 0 & 0 & 0 & 0 & 0 & 0 & 0 & 0 & 0 & 0 & 0 \\
t_{0,34} & t_{1,34} & t_{2,34} & t_{12,34} & -t_{0,4} & -t_{1,4} & -t_{2,4} & -t_{12,4} & t_{0,3} & t_{1 ,3} & t_{2,3} & t_{12,3} & 0 & t_{1,0}  & t_{2,0} & t_ {12,0}\\
0 & t_{0,34} & 0 & -t_{2,34} & 0 & t_{0,4} & 0 & -t_{2,4} & 0 & -t_{0,3}  & 0 & t_{2,3} & 0 & 0 & 0 & -t_{2,0}\\
0 & 0 & t_{0,34} & t_{1,34} & 0 & 0 & t_{0,4} & t_{1,4} & 0 & 0 & -t_{0,3} & -t_{1,3} & 0 & 0 & 0 & t_{1,0} \\
0 & 0 & 0 & t_{0,34}  & 0 & 0 & 0 & -t_{0,4} & 0 & 0 & 0 & t_{0,3} & 0 & 0 & 0 & 0
\end{array}
\right].
$$ 
\end{tiny}
The level zero variable $t_{0,0}$ does not appear and the level four variable $t_{12,34}$ only appears once. The matrix is not lower triangular, unlike the cluster matrices defined in \cite[Section~2]{FSS}. However, it is nilpotent of order $2$, so $T(t)^3 = 0$. The fourth column of the exponential matrix is of the form
\begin{footnotesize}
$$
\psi = \exp(T(t))e_{12} = 
    \begin{bmatrix}
         t_{12,0}\\
-t_{2,0}\\
t_{1,0}\\
1\\
t_{0,3}t_{12,0}-t_{2,0}t_{1,3}+t_{1,0}t_{2,3}+t_{12,3}\\
t_{2,3}\\
-t_{1,3}\\
t_{0,3}\\
t_{0,4}t_{12,0}-t_{2,0}t_{1,4}+t_{1,0}t_{2,4}+t_{12,4}\\
t_{2,4}\\
-t_{1,4}\\
t_{0,4}\\
t_{12,0}t_{0,34}-t_{1,4}t_{2,3}+t_{1,3}t_{2,4}+t_{12,34}\\
-t_{2,0}t_{0,34}+t_{0,4}t_{2,3}-t_{0,3}t_{2,4}-t_{2,34}\\
t_{1,0}t_{0,34}-t_{0,4}t_{1,3}+t_{0,3}t_{1,4}+t_{1,34}\\
t_{0,34}
    \end{bmatrix}.
$$  
\end{footnotesize}
\end{example}

\begin{proposition}
The exponential parametrization is bijective and has a polynomial inverse.
\end{proposition}

\begin{proof}
First we notice that at level zero we have $t_{0, 0} = \psi_{[d]} = 1$ and at level one we have $t_{j, 0} = \pm \psi_{[d] \backslash \{j\}}$ and $t_{0, b} = \pm \psi_{[d] \cup \{b\}}$. If  $\psi_{J}$ has level $r$ then we can write $\psi_J$ as $\pm t_{I,B}$, where $I = [d] \backslash J$ and $B = J \backslash [d]$, plus a polynomial in variables $t$ of level $< r$. Each of these lower level $t$'s can now be replaced with a polynomial in $\psi$, by the induction hypothesis. This yields a representation for $t_{I,B}$ as $\pm \psi_J$ plus a polynomial in lower level $\psi$-coordinates.
\end{proof}

The polynomials $\psi_{[2d] \backslash [d]}$ are called the \textit{master polynomials}. For example, when $d = 3$ and $n = 6$, we have the following master polynomial $\psi_{456}(t)$:
\begin{align*}
     t_{23, 0}t_{0, 56}t_{1, 4}-t_{23, 0}t_{0, 46}t_{1, 5}+t_{23, 0}t_{0, 45}t_{1, 6}
-t_{13, 0}t_{0, 56}t_{2, 4}+t_{13, 0}t_{0, 46}t_{2, 5}-t_{13, 0}t_{0, 45}t_{2, 6}\\
+t_{12, 0}t_{0, 56}t_{3, 4}
-t_{1, 6}t_{2, 5}t_{3, 4}+t_{1, 5}t_{2, 6}t_{3, 4}
-t_{12, 0}t_{0, 46}t_{3, 5}+t_{1, 6}t_{2, 4}t_{3, 5}-t_{1, 4}t_{2, 6}t_{3, 5}\\
+t_{12, 0}t_{0, 45}t_{3, 6}-t_{1, 5}t_{2, 4}t_{3, 6}
+t_{1, 4}t_{2, 5}t_{3, 6}
+t_{0, 56}t_{123, 4}-t_{0, 46}t_{123, 5}+t_{0, 45}t_{123, 6}\\+t_{23, 0}t_{1, 456}
-t_{13, 0}t_{2, 456}+t_{12, 0}t_{3, 456}
+t_{3, 6}t_{12, 45}-t_{3, 5}t_{12, 46}
+t_{3, 4}t_{12, 56}\\-t_{2, 6}t_{13, 45}+t_{2, 5}t_{13, 46}-t_{2, 4}t_{13, 56}
+t_{1, 6}t_{23, 45}
-t_{1, 5}t_{23, 46}+t_{1, 4}t_{23, 56}+t_{123, 456}.
\end{align*}
It has $31$ terms. We compare that number to the $16$ terms the master polynomial defined in \cite[Section~2]{FSS} has. All the terms of the master polynomial in first quantization do appear in the master polynomial above. However, in second quantization we have extra terms, like $t_{23, 0}t_{0, 56}t_{1, 4}$, which have variables indexed by sets $I$ and $B$ where $|I| \ne |B|$. The master polynomial in first quantization is indexed by all uniform block permutations of $[2d]$.
The variables $t$ in the master polynomials above are also indexed by elements in the set $[2d]$. More specifically, the monomials of the master polynomial are indexed by set partitions of $[2d]$.

The other polynomials $\psi_I$ in the exponential parametrization are just relabelings of the master polynomials. The variables $t$ in $\psi_I(t)$ are indexed by elements in the symmetric difference $I \oplus [d]$, and the monomials are indexed by set partitions of $I \oplus [d]$. In particular, if $|I \oplus [d]| = 2k$ is even, $\psi_I$ is a relabeling of master polynomial $\psi_{[2k] \backslash [k]}$. If $|I \oplus [d]| = 2k - 1$ is odd, we add a formal element $\varnothing$ to the index set, and the polynomial $\psi_I$ also becomes a relabeling of the master polynomial $\psi_{[2k] \backslash [k]}$.

\begin{example}[$n = 5$, $d = 2$]
There are $32$ polynomials $\psi_I$ in the exponential parametrization. Polynomials like

\begin{footnotesize}
\begin{align*}
    &\psi_{35}(t) = t_{12,0}t_{0,35} -t_{1,5}t_{2,3}+t_{1,3}t_{2,5}+t_{12,35}, \,\,\,\, \psi_{1345}(t) = t_{0,45}t_{2,3}-t_{0,35}t_{2,4}+t_{0,34}t_{2,5}+t_{2,345},\\
     &\psi_{45}(t) = t_{12,0}t_{0,45} -t_{1,5}t_{2,4}+t_{1,4}t_{2,5}+t_{12,45}, \,\,\,\,
    \psi_{2345}(t) = t_{0,45}t_{1,3}-t_{0,35}t_{1,4}+t_{0,34}t_{1,5}+t_{1,345}
\end{align*}
\end{footnotesize}
\vspace{-0em}

\noindent are relabelings of the master polynomial $\psi_{34}(t) = t_{12,0}t_{0,34}-t_{1,4}t_{2,3}+t_{1,3}t_{2,4}+t_{12,34}$. Following polynomials are also relabelings of $\psi_{34}(t)$, but a formal element is added to the index set to create a one-to-one correspondence with $\{1,2,3,4\}$:

\begin{small}
    \begin{align*}
    &\psi_5(t) = t_{12,0}t_{0,5}-t_{1,5}t_{2,0}+t_{1,0}t_{2,5}+t_{12,5}, \,\,\,
     \psi_{235}(t) = t_{1,0}t_{0,35} -t_{0,5}t_{1,3}+t_{0,3}t_{1,5}+t_{1,35},\\
      &\psi_{12345}(t) = t_{0,34}t_{0,5}-t_{0,35}t_{0,4} +t_{0,45}t_{0,3}+t_{0,345}.
\end{align*}
\end{small}
\vspace{-0.5em}

\noindent The formal element is not shown when writing out the polynomials.
The polynomial $\psi_{345}(t)$ has $31$ terms and is a relabeling of the master polynomial $\psi_{456}(t)$.
    
\end{example}

By looking at the number of terms in the master polynomials we get the sequence
\begin{align*}
    1, 4, 31, 379, 6\,556, 150\,349, 4\,373\,461, \dots \tag{A005046}
\end{align*}
The sequence enumerates the partitions of an even set into even blocks. This suggests we can write polynomials $\psi_I$ as a linear combination of monomials indexed by even set partitions. 
Write monomials corresponding to set partition $\pi = (\pi_1, \dots, \pi_k)$ as
\begin{equation}
    t_\pi = t_{\pi_1 \cap [d], \pi_1 \backslash [d]}\cdots t_{\pi_k \cap [d], \pi_k \backslash [d]} \,\,\, \text{or}\,\,\, x_{\pi} = x_{\pi_1 \oplus [d]} \cdots x_{\pi_k \oplus [d]}.
\end{equation}

\begin{theorem}\label{thm:psiparam}
    The master polynomials $\psi_{[2d] \backslash [d]}$ are linear combinations of the monomials $t_\pi$ where $\pi$ is an even partition of $[2d]$. More specifically, we get
    $$
    \psi_{[2d] \backslash [d]}(t) = \sum_{\substack{ \pi \, \vdash \, [2d]\\ \pi \text{ even}}} \operatorname{sign}(\pi) t_\pi
    $$
    where the sign of $\pi$ is the product of the signs of the following two permutations
    $$
    [d] \mapsto (\pi_1 \cap [d], \dots, \pi_k \cap [d]), \quad [2d] \backslash [d] \mapsto (\pi_1 \backslash [d], \dots, \pi_k \backslash [d]).
    $$
\end{theorem}

\begin{proof}
    By Remark \ref{re:marepFD} we can write
    \begin{align*}
        &\psi_{[2d] \backslash [d]}(t) = e_{[2d] \backslash [d]}^\dagger \exp(T(t))e_{[d]}\\
        &= \sum_{k = 1}^d \frac{1}{k!}\sum_{\substack{I_i \subseteq [d],\\ B_i \subseteq [2d] \backslash [d]}} t_{I_1, B_1} \cdots t_{I_k, B_k} 1^\dagger a_{2d} \cdots a_{d + 1} a_{B_1}^\dagger a_{I_1} \cdots  a_{B_k}^\dagger a_{I_k} a_1^\dagger \cdots a_d^\dagger 1.
    \end{align*}
    By Remark \ref{re:constSR}, we only need to sum over index sets $\{(I_i, B_i) \, : \, 1 \le i \le k \}$ where the disjoint union $\sqcup_i (I_i \sqcup B_i)$ is equal to $[2d]$.
    We also note that if there are two pairs $(I, B)$ and $(I', B')$ where $|I| + |B|$ and $|I'| + |B'|$ are odd, then $a_{B}^\dagger a_{I}$ and $a_{B'}^\dagger a_{I'}$ anti commute. If they are both even they commute. Thus we are summing over the even set partitions of $[2d]$, that is
    $$
    \psi_{[2d] \backslash [d]} = \sum_{\substack{ \pi \, \vdash \, [2d]\\ \pi \text{ even}}} t_\pi 1^\dagger a_{2d} \cdots a_{d + 1} a_{B_1}^\dagger a_{I_1} \cdots  a_{B_k}^\dagger a_{I_k} a_1^\dagger \cdots a_d^\dagger 1.
    $$
    The number of crossings in $a_{2d} \cdots a_{d + 1} a_{B_1}^\dagger a_{I_1} \cdots  a_{B_k}^\dagger a_{I_k} a_1^\dagger \cdots a_d^\dagger$ is the sum of the number of inversions of the two permutations defined in the statement of the theorem.
\end{proof}

As a direct result, we see that the degree of the polynomial $\psi_I$ depends only on the symmetric difference $I \oplus [d]$. That is $\deg \psi_I = \lceil \frac{1}{2} |I \oplus [d]| \rceil$.
The exponential parametrization is a nonlinear bijection and its inverse is also defined by polynomials. Just as for $\psi_I$ we see that each polynomial $t_{I,B}(\psi) = x_J(\psi)$ is just a relabeling of the master polynomial $x_{[2d] \backslash [d]}$. We use the coordinates $c_{I,B}$ instead of $\psi_J$ and write monomials for $c$ corresponding to partition $\pi = (\pi_1, \dots, \pi_k)$ as we do for $t$. Then for each even partition $\rho$ we get that
$$
(-1)^\nu \operatorname{sign}(\rho)c_\rho = \sum_{\pi \le \rho} \operatorname{sign}(\pi)t_\pi
$$
where $\nu = \binom{d}{2} - \sum_{r = 1}^k \binom{|\rho_r \cap d|}{2}$. See the proof of \cite[Theorem~2.5]{FSS} where the sign $\nu$ appears as well. The Möbius function for even partitions of $[2d]$ is $\mu(\pi) = (-1)^{k - 1}(k - 1)!$ so using Möbius inversion we obtain
\begin{equation}
    x_{[2d] \backslash [d]}(c) = \sum_{\substack{ \pi \, \vdash \, [2d]\\ \pi \text{ even}}} (-1)^{\nu + k - 1}(k  - 1)!\operatorname{sign}(\pi) c_\pi.
\end{equation}
The proof of this statement is analogous to the proof of \cite[Theorem 2.5]{FSS}.

\section{Fock space truncation varieties}\label{se:TruncVars}

We define the \textit{truncation grid}
$$
\mathcal{G}
= 
\{
(m,\ell) ~:~
0 \le m \le d \quad {\rm and} \quad 0 \le \ell \le n - d
\} \backslash \{(0,0)\}
$$
enumerating feasible annihilation and creation levels. The corresponding construct in first quantization is the diagonal of the truncation grid or  the set $[d] = \{1,2,\dots,d\}$.
For a proper subset $\sigma  \subsetneq \mathcal{G}$, called the \textit{level set}, we define a subspace $\mathcal{V}_\sigma$ of $\mathcal{V}$ spanned by occupational number states $e_J$ with excitation levels $(a(J) , c(J))$ in $\sigma$. In symbols:
$$
\mathcal{V}_\sigma = \operatorname{span}\{\, e_J \, : J \subseteq [n] \,\text{ and }\, \, (a(J), c(J)) \,\in\, \sigma \,\}.
$$
The projection of cluster amplitude $t \in \mathcal{V}$ onto $\mathcal{V}_\sigma$ is denoted $t_\sigma$. As a variety in the projectivization of $\mathcal{V}$, $\PP(\mathcal{V}) \cong \PP^{2^n - 1}$, it is the vanishing set of the linear ideal
$$
\langle \, x_J : J \subseteq [n] \,\text{ and }\, (a(J) , c(J))\, \in \,\mathcal{G} \backslash \sigma \,\rangle.
$$
We look at the restriction of the exponential map to the subspace $\mathcal{V}_\sigma$:
$$
\mathcal{V}_\sigma \to \mathcal{F}, \quad t_\sigma \mapsto \psi = \exp(T(t_\sigma))e_{[d]}.
$$
It is injective and maps $\mathcal{V}_\sigma$ into the Fock space $\mathcal{F}$, which further maps to the projective space $\PP(\mathcal{F}) \cong \PP^{2^n - 1}$ using homogenizing variable $\psi_{[d]}$.
The \textit{Fock space truncation variety} $V_\sigma$ is defined as the closure of the image of $\mathcal{V}_\sigma$ under this map to $\PP^{2^n - 1}$. The varieties live in the projective space of binary tensors and contain the truncated quantum states. 
Since the exponential map is injective, the dimension of the truncation variety is equal to the dimension of the subspace $\mathcal{V}_\sigma$. In symbols:
$$
\dim(V_\sigma) = \dim(\mathcal{V}_\sigma) = | \{ J : (a(J), c(J)) \in \sigma\}| = \sum_{(m,\ell) \in \sigma} \binom{d}{m} \binom{n - d}{\ell}.
$$

In the case when the level set $\sigma$ is a subset of the diagonal of the truncation grid, i.e. $\sigma = \{ (k,k) \, :\, k \}$ we obtain the truncation varieties from \cite{FSS}, embedded in the projective space of binary tensors $\PP^{2^n - 1}$. Those truncation varieties contain quantum states of systems with a fixed number of electrons. When we relax the particle-number conservation we work with truncation varieties corresponding to level sets not lying on the diagonal. In FSCC we truncate ionized quantum states using for example CCSD.  These quantum states live in the truncation variety corresponding to $\sigma = \{(1,0), (2,1), (1,1), (2,2)\}$. The truncation varieties defined in this paper extend the truncation varieties defined in \cite[Section~3]{FSS} from fixed-N CC to FSCC.

\begin{example}
    For $d = 2$ and $n = 4$, the truncation grid is equal to
    $$
    \mathcal{G} = \{(0,1), (0,2), (1,0), (1,1), (1,2), (2,0), (2,1), (2,2)\}.
    $$
    This is a set of size $8$. Therefore, there are $2^8 - 2 = 254$ choices of non-empty proper level sets $\sigma \subsetneq \mathcal{G}$. Hence, we have $254$ truncation varieties in $\PP^{15}$ corresponding to $2$ base electrons in $4$ orbitals. The truncation varieties corresponding to $\sigma = \{(1,1)\}, \{(2,2)\}, \{(1,1), (2,2) \}$ are respectively isomorphic to the truncation varieties $V_{\{1\}}, V_{\{2\}}$ and $V_{\{1,2\}}$ defined in \cite[Section~3]{FSS} embedded into $\PP^{15}$. We do not obtain $254$ distinct truncation varieties since some of the varieties are isomorphic. For example, the truncation variety corresponding to $\{(0,1), (1,1)\}$ is isomorphic to the truncation variety corresponding to $\{(1,0), (1,1)\}$. We also note that $119$ of the $254$ truncation varieties are linear subspaces, see upcoming Theorem~\ref{thm:linear}. Also, $74$ truncation varieties fulfill the hypothesis of Theorem \ref{thm:degofgraph} and their CC degree is the degree of the graph of the exponential map.
\end{example}

For $d$ base electrons and $n$ orbitals, the truncation grid $\mathcal{G}$ has cardinality $(d + 1)(n - d + 1) - 1 = dn + n - d^2$. The number of truncation varieties $V_\sigma$ is $2^{dn + n - d^2} - 2$. Some of these varieties are isomorphic.
With that being said, we present an isomorphism of truncation varieties, known as the particle-hole formalism:

\begin{proposition}\label{prop:parthole}
    Fix a proper subset $\sigma \subsetneq \mathcal{G}$ and let $n \ge 2d$. There is a linear isomorphism between truncation variety $V_\sigma$ for $(d,n)$ and $V_{\overline{\sigma}}$ for $(n - d, n)$. Here
    $
    \overline{\sigma} = \{ (\ell, m) \, : \, (m, \ell) \in \sigma \}.
    $
\end{proposition}

\begin{proof}
    This isomorphism is obtained through the following relabeling of Plücker coordinates:
    $$
    J \mapsto J' = \{n + 1 - j : j \notin J \}.
    $$
    This is the same relabeling used in an analogous isomorphism in \cite[Proposition~3.7]{FSS}.
\end{proof}
Let $\overline{x}_I(\psi)$ denote the homogenization of $x_I(\psi)$ using homogenizing variable $\psi_{[d]}$.

\begin{theorem}
    The homogeneous prime ideal of the truncation variety $V_\sigma \subseteq \PP^{2^n - 1}$ is the following saturation \cite[Section~4.4]{cox1997ideals}:
    $$
    I(V_\sigma)\, =\, \langle \, \overline{x}_I(\psi) \, : \, (a(I), c(I)) \in \mathcal{G} \backslash \sigma \, \rangle : \langle\, \psi_{[d]} \,\rangle^\infty.
    $$
    In particular, the ideal of the restriction of $\,V_\sigma$ to the affine chart $\,\CC^{2^n -1} = \{ \psi_{[d]} =1 \}$ of projective space $\,\PP^{2^n-1}$ is the complete intersection 
$$
I(V_\sigma) + \langle \, \psi_{[d]} - 1 \, \rangle = \langle \, x_I(\psi) \,:\,  (a(I), c(I)) \in \mathcal{G} \backslash \sigma \,\rangle.
$$
\end{theorem}

For the sake of the next theorem, we define a partial order $\preceq$ on the truncation grid $\mathcal{G}$.
We say excitation level $(m, \ell)$ is less than excitation level $(k, p)$, $(m, \ell) \preceq (k, p)$ if the creation and annihilation levels fulfill $m \le k$, $\ell \le p$ and either the level $m + \ell$ is even or level $k + p$ is odd. One can check that this is in fact a partial order. The level sets $\sigma$ and $\mathcal{G} \backslash \sigma$ can be defined as posets with the same partial order $\preceq$.
We say $\sigma \subseteq \mathcal{G}$ is \textit{closed under addition with respect to partial order $\preceq$} if for $(m, \ell), (k, p) \in \sigma$, where $(m, \ell), (k, p) \preceq (m + k, \ell + p) \in \mathcal{G}$, then $(m + k, \ell + p) \in \sigma$. The following result generalizes \cite[Theorem~3.10]{FSS} where we assume particle-number conservation. There $\sigma$ is a totally ordered set in a totally ordered truncation grid $[d]$. We emphasize that the truncation grid should not be confused with the posets introduced in \cite{vcivzek1966correlation}. In that setting, a partial order is defined on the indexing sets $I \in \binom{[n]}{d}$, whereas here we define a partial order on the levels of the indexing set $I \subseteq [n]$.

\begin{theorem}\label{thm:linear}
    The variety $V_\sigma$ is linear if and only if $\sigma$ is closed under addition with respect to the partial order defined above.
    In other words, $V_\sigma$ is linear if and only if for all $(m,\ell), (k,p) \in \sigma$ such that either $m + \ell$ or $k + p$ is even, then either
    $$
    (m, \ell) + (k, p) \notin \mathcal{G} \,\,\text{ or }\,\,(m,\ell) + (k,p) \in \sigma.
    $$
\end{theorem}

\begin{proof}
    We identify $V_\sigma$ with its restriction to  the affine chart $\mathcal{F}' = \CC^{2^n-1}$. We assume that $\sigma$ is closed under addition with respect to $\preceq$.
    We first prove by induction over poset $\mathcal{G} \backslash \sigma$ that $\psi_I = 0$ for all $I \subseteq [n]$ with creation and annihilation level in $\mathcal{G} \backslash \sigma$.

    Since $\sigma$ is closed under addition with respect to $\preceq$, the minimal elements of $\mathcal{G} \backslash \sigma$ are a subset of the minimal elements of $\mathcal{G}$. The minimal elements of $\mathcal{G}$ are $(0,1), (1,0), (2,0), (0,2)$ and $(1,1)$. These are also all pairs of level $\le 2$.
    We take a set $I$ of creation and annihilation level $(m, \ell)$ that is minimal in $\mathcal{G} \backslash \sigma$. Then $x_I(\psi) = \psi_I = 0$, since $I$ is of level $\le 2$.

    Now consider a set $K$ with creation and annihilation level $(m, \ell) \in \mathcal{G} \backslash \sigma$. Then
    $$
    x_K(\psi) = \sum_{\substack{ \pi \, \vdash \, K \oplus [d]\\ \pi \text{ even}}} \pm \psi_\pi =  \pm \psi_K + \sum_{j} a_j \psi_{K_1^{(j)}} \cdots \psi_{K_{r_j}^{(j)} } = 0
    $$
    where $a_j \in \ZZ^\ast$ and $K_s^{(j)} = \pi_s^{(j)} \oplus [d]$ where $\{\pi_1^{(j)}, \dots, \pi_{r_j}^{(j)} \}$ forms an even set partition of $K \oplus [d]$.  We let $(m_s^{(j)}, \ell_s^{(j)})$ be the creation and annihilation levels of the sets $K_{s}^{(j)}$, and notice that $\sum_{s = 1}^{r_j} (m_s^{(j)}, \ell_s^{(j)}) = (m,\ell)$ for each $j$.
    If $K$ has an odd level then $(m_s^{(j)}, \ell_s^{(j)}) \preceq (m, \ell)$ for all $s$ and $j$. 
    If $K$ has an even level, then all the $(m_s^{(j)},\ell_s^{(j)})$ have an even level, so $(m_s^{(j)}, \ell_s^{(j)}) \preceq (m, \ell)$ as well. Since $\sigma$ is closed under addition with respect to $\preceq$, there must be one pair $(m_s^{(j)}, \ell_s^{(j)}) \in \mathcal{G} \backslash \sigma$ for each $j$. Thus by induction there is some $s$ for each $j$ such that $\psi_{K_s^{(j)}} = 0$ and hence $0 = x_K(\psi) = \psi_K$. The defining equations of the truncation variety are therefore linear of the form $x_K(\psi) = \psi_K$ where $K \in \mathcal{G} \backslash \sigma$.

    We now assume $\sigma$ is not closed under addition with respect to $\preceq$. 
    We take a minimal $(m, \ell) \in \mathcal{G} \backslash \sigma$, such that there are some $(i,r), (j,s) \in \sigma$ where $(i,r), (j,s) \preceq (m, \ell)$ in $\mathcal{G}$ and $(m, \ell) = (i,r) + (j,s)$.
    Consider a polynomial $x_K(\psi)$ vanishing on $V_\sigma$ where $K$ has creation and annihilation level $(m,\ell)$. Fix any degree-compatible monomial order. The initial monomial of $x_K(\psi)$ has degree $> 1$,
    $$
    \operatorname{in}(x_K(\psi)) = \psi_{K_1} \cdots \psi_{K_{r}}, \quad \text{ where } r \ge 2. 
    $$
    Assume some element in the initial ideal of $\mathcal{I}(V_\sigma)$ divides $\operatorname{in}(x_K(\psi)) = \psi_{K_1} \cdots \psi_{K_{r}}$. This element must divide a monomial of some generator $x_I(\psi)$; here $(a(I), c(I)) \in \mathcal{G} \backslash \sigma$. Without loss of generality we may pick the divisor to be a monomial of $x_I(\psi)$. We also see that $(a(I), c(I))  \preceq  (m,\ell)$, so by the same argument as above we get $x_I(\psi) = \psi_I$. We also obtain the relation $\psi_I = 0$ and thus $\psi_I \nmid \operatorname{in}(x_K(\psi))$, a contradiction. Hence $\operatorname{in}(x_K(\psi))$ is a minimal generator for the initial ideal of $\,\mathcal{I}(V_\sigma)$. This cannot be an initial ideal for a linear variety, and therefore $V_\sigma$ itself is not linear.
\end{proof}

\begin{example}[$d = 3$, $n = 6$] 
There are $2^{15} - 2 = 32\,766$ truncation varieties in $\PP^{63}$. Among those, $4\,790$ are linear by Theorem \ref{thm:linear} and $2 \,186$ fulfill the hypothesis of Theorem \ref{thm:degofgraph}.  We look at a few chosen cases of non-linear truncation varieties:

\begin{enumerate}
    \item[]$\sigma = \{(1,0), (2,0)\}$:\\
    This is a hypersurface of the subspace generated by coordinates $\psi_I$ where $I \subseteq [3]$. Its defining equation is the homogeneous inverse coordinate $\overline{x}_{\varnothing}(\psi)$.\\
    \item[]$\sigma = \{ (1,0), (2,0), (1,1),(0,1), (0,2)\}$.\\
    Its ideal is generated by $364$ quadrics. It is of dimension $21$ and degree $33\,592$.\\
    \item[]$\sigma = \{(1,0), (2,1),(1,1), (2,2)\}$.\\
    Its ideal is generated by $7$ quadrics and the $30$ coordinates indexed by sets of size $0,1,4,5$ or~$6$. It is of dimension $30$ and degree $43$, see Example 7.3.\\
    \item[] $\sigma = \{(1,0), (1,1), (0,1)\}:$\\
    This is the partial flag variety $\mathcal{F}\ell (2,3,4,6)$. It has dimension $15$ and degree $4\,550$. It is the zero set of $281$ quadrics and the $14$ coordinates indexed by sets of size $0,1,5$ or~$6$.\\
    \item[] $\sigma = \{(2,0), (1,1), (0,2)\}$:\\
    This is the spinor variety $S_+$ for a $12$ dimensional space like $\CC^{12}$. It has dimension $15$ and degree $286$. It is the zero set of $66$ quadrics and the $32$ coordinates indexed by even sets.    
\end{enumerate}
    
\end{example}

When extending from the exterior power $\mathcal{H} \cong \wedge^d \mathbb{R}^n$ to the full Fock space $\mathcal{F}$, we allow the number of particles in the electronic system to vary. In this broader setting, we encounter several well-known algebraic varieties, such as the flag variety and the spinor variety, which will be examined in detail in the next section. Beyond these classical examples, the framework developed here naturally gives rise to many new and intriguing varieties, whose study is of independent mathematical interest, regardless of their origins in quantum chemistry. At the same time, returning to the motivating context, it is natural to ask which of these Fock space truncation varieties are most relevant for practical applications in quantum chemistry.

\begin{remark}[EOM-CC for ionization and electron attachment]\label{re:eomcc}
In addition to FSCC, the equation-of-motion coupled cluster method (EOM-CC) also addresses ionized electronic systems and electron attachment; see \cite[Chapter~13.4]{shavitt2009many}. There we introduce the \textit{ionization operator}:
\begin{equation}
    T(t_\iota) = \sum_{\substack{ I \subseteq [d], \,  B \subseteq [n] \backslash [d] \\ |I| = |B| + 1}} t_{I,B} a_{b_\ell}^\dagger...a_{b_1}^\dagger
a_{i_1}...a_{i_m} =  \sum_{\substack{ I \subseteq [d], \,  B \subseteq [n] \backslash [d] \\  |I| = |B| + 1}} t_{I,B} a_B^\dagger
a_I
\end{equation}
where $\iota = \{(k + 1, k), 0 \le k \le d \}$, is the sub diagonal of the truncation grid.
This operator is an element of the Fermi-Dirac algebra. 
The ionized quantum states are then parameterized by the following exponential map
\begin{equation}\label{eq:ionparam}
    \mathcal{V_\iota} \times \mathcal{V}_{[d]} \to \wedge^{d - 1} \RR^n, \quad (t_\rho, t_{[d]}) \mapsto \psi = T(t_\iota) \exp{T(t_{[d]})}e_{[d]}.
\end{equation}
Here $[d] = \{(k,k) : 1 \le k \le d\}$ is the diagonal of the truncation grid. 
We truncate the ionized quantum states using subsets of the sub diagonal $\iota$ and the diagonal $[d]$ and restricting the exponential map (\ref{eq:ionparam}) accordingly. The closure of the image of this restriction is a variety of $\wedge^{d - 1} \mathbb{R}^n$. In most applications, we fix a positive integer $k$ and truncate with subsets $\{(i + 1, i) : 0 \le i \le k \} \subseteq \iota$ and $[k] \subseteq [d]$.
Electron attachment is treated dually using the particle-hole formalism. We define the \textit{electron attachment operator} as $T(t_\epsilon)$ where $\epsilon = \overline{\iota} = \{(k, k + 1), 0 \le k \le d \}$ is the super diagonal.
\end{remark}

At first glance, the varieties in Remark \ref{re:eomcc} arising from EOM-CC appear incompatible with our construction of truncation varieties. The following theorem, however, demonstrates how these two frameworks are related.

\begin{theorem}\label{thm:eomcc}
   Let $\tau \subseteq \iota$ be a set on the sub diagonal and $\sigma \subseteq [d]$ a set on the diagonal of the truncation grid. The variety of truncated ionized EOM-CC quantum states, parameterized by the restriction of (\ref{eq:ionparam}) to $\mathcal{V}_\tau \times \mathcal{V}_\sigma$ is equal to the truncation variety $V_{\tau \cup \sigma}$ restricted to the $(d - 1)$st exterior power $\wedge^{d - 1}\RR^n$. A dual statement holds for electron attachment.
\end{theorem}

\begin{proof}
    The truncation variety $V_{\tau \cup \sigma}$ is parameterized by the following restriction of the exponential map:
    $$
    \mathcal{V}_{\tau \cup \sigma} \to \mathcal{F}, \quad  t_{\tau \cup \sigma} \mapsto \psi = \exp{(T(t_\tau) + T(t_{\sigma}))}e_{[d]}.
    $$
    The operator $T(t_{\sigma})$ is a linear combination of monomials of even degree, hence we see by Remark \ref{re:exteriorsofFD} that $T(t_{\sigma})$ and $T(t_{\tau})$ commute. Thus we can factor the exponential:
    $$
    \exp{(T(t_\tau) + T(t_{\sigma}))} = \exp{T(t_\tau)}\exp{T(t_\sigma)}.
    $$ 
    Moreover, the operator $T(t_{\tau})$ is a linear combination of monomials of odd degree so, by Remark \ref{re:exteriorsofFD},  $T(t_{\tau})$ is nilpotent of order $1$. We can write $\exp{T(t_\tau)} = 1 + T(t_\tau)$. We rewrite the restricted exponential map:
    $$
    \mathcal{V}_{\tau \cup \sigma} \to \mathcal{F}, \quad  t_{\tau \cup \sigma} \mapsto \psi =  \exp{T(t_{\sigma})}e_{[d]} + T(t_{\tau})\exp{T(t_{\sigma})}e_{[d]}.
    $$
    The first summand of the parametrization is an element in the $d$th exterior power $\wedge^d \RR^n$ and the second summand is an element of $\wedge^{d - 1} \RR^n$ proving our statement.
\end{proof}

\section{Flag varieties and spinor varieties}\label{se:FlagsSpinors}

We define a \textit{flag} in the vector space $\CC^n$ to be a sequence of subspaces increasing with respect to inclusion, that is: 
$$
F_\bullet: \{0\} = F_0 \subsetneq F_1 \subsetneq \cdots \subsetneq F_k = \CC^n.
$$
We say a flag is \textit{complete} if $d_i = \dim (F_i) = i$ for all $0 \le i \le k$; otherwise, it is called a \textit{partial} flag. The \textit{partial flag variety } $\mathcal{F}l(d_1, \dots, d_k, n) \subseteq \PP^{2^n - 1}$ is the set of all flags in $\CC^n$ with dimensions $d_1 \le d_2 \le \cdots \le d_k$. 

We recall that a $d$ dimensional subspace of $\CC^n$ can be parametrized by the maximal minors of a $d \times n$ matrix $\Theta$ with entries in $\CC$. The matrix $\Theta$ is of full rank and since row reduction does not change the row span of $\Theta$ we may assume that
$
\Theta  = [I_d \, | \, M]
$
where $M$ is a $d \times (n - d)$ matrix. 
Similarly, for a flag $F_\bullet$ we look at a full rank $d_k \times n$ matrix $\Theta = [I_{d_1} | M]$ where the first $d_i$ rows of $\Theta$ span $F_i$. The flag $F_\bullet$ can then be parametrized by the maximal minors of the first $d_i$ rows of $\Theta$ for each $i$. We refer to Ezra Miller's and Bernd Sturmfels' book \cite[Section~14.1]{miller2005combinatorial} for a more detailed introduction to flag varieties.

\begin{remark}\label{re:paramflag}
    We describe an alternative parametrization of the flag variety $\mathcal{F}l(d - 1, d, d + 1, n)$. We look at the $(d + 1) \times (n + 1)$ matrix 
    $$
    \Theta = \begin{bmatrix}
        1 & 0 & \cdots & 0 & \theta_{1,d + 1} & \cdots & \theta_{1,n} & \theta_{1,n + 1}\\
        0 & 1 & \cdots & 0 & \theta_{2,d + 1} & \cdots & \theta_{2,n} & \theta_{2,n + 1}\\
        \vdots & \vdots & \ddots & \vdots & \vdots & \ddots & \vdots & \vdots\\
        0 & 0 & \cdots & 1 & \theta_{d,d + 1} & \cdots & \theta_{d,n} & \theta_{d,n + 1}\\
        0 & 0 & \cdots & 0 & \theta_{d + 1,d + 1} & \cdots & \theta_{d + 1,n} & 0\\
    \end{bmatrix}.
    $$
    The first $d$ rows of $\Theta$ span a $d$ dimensional subspace $E$ of $\CC^{n + 1}$. The intersection of $E$ with the orthogonal space $(e_{n + 1})^\perp \cong \CC^n$ is generically isomorphic to a $d - 1$ dimensional subspace of $\CC^n$, which we will denote by $F_{d - 1}$. 
    We now look at a submatrix $\Tilde{\Theta}$ of $\Theta$ taking only the first $n$ columns. Its row span is a subspace $F_{d + 1}$ of $\CC^n$, which is generically of dimension $d + 1$. Also, the first $d$ rows of $\Tilde{\Theta}$ span a $d$ dimensional subspace $F_{d} \subseteq F_{d + 1}$. Therefore, we have obtained a partial flag 
    $$
    F_{d - 1} \subseteq F_d \subseteq F_{d + 1} \subseteq \CC^n.
    $$
   This flag is parametrized by the $d \times d$ minors of the first $d$ rows of $\Theta$ and the $(d + 1) \times (d + 1)$ minors of the first $n$ columns of $\Theta$. These minors define a birational parametrization of the partial flag variety $\mathcal{F}l(d - 1, d, d + 1, n)$.
\end{remark}

The partial flag variety $\mathcal{F}l(d - 1, d, d + 1,n)$ appears as a truncation variety. Before stating the general theorem, we first show this in an example
  
\begin{example}[A small flag]
We look at the case when $d = 2$, $n = 4$ and $\sigma = \{(1,1), (0,1), (1,0)\}$. The cluster matrix $T(t_\sigma)$ is a $16 \times 16$ matrix, nilpotent of order $2$.
The truncation variety is of dimension
$$
\dim(V_\sigma) = \binom{2}{1} + \binom{2}{1} + \binom{2}{1} \cdot \binom{2}{1}  = 8
$$
and degree $\deg(V_\sigma) = 12$. The coordinates of the parametrization are
\begin{align*}
    &\psi_\varnothing = 0, \, \psi_{1} = -t_{2,0}, \, \psi_{2} = t_{1,0}, \, \psi_{12} = 1, \, \psi_{3} = -t_{2,0}t_{1,3}+t_{1,0}t_{2,3}, \,  \psi_{13} = t_{2,3},\\ &\psi_{2,3} = -t_{1,3}, \, \psi_{123} = t_{0,3}, \, \psi_{4} = -t_{2,0}t_{1,4}+t_{1,0}t_{2,4}, \, \psi_{14} = t_{2,4}, \, \psi_{24} = -t_{1,4},\\ &\psi_{124} = t_{0,4}, \, \psi_{34} = -t_{1,4}t_{2,3}+t_{1,3}t_{2,4}, \, \psi_{134} = -t_{2,4}t_{0,3}+t_{2,3}t_{0,4}, \\ &\psi_{234} = t_{1,4}t_{0,3}-t_{1,3}t_{0,4}, \, \psi_{1234} = 0.
\end{align*}
These are the $2 \times 2$ minors of the first two rows; and the $3 \times 3$ minors of the last four columns of the $3 \times 5$ matrix
\begin{equation}\label{eq:flag24param}
     M = \begin{bmatrix}
         0 & 0 & 0  & t_{0,3} & t_{0,4}\\
        t_{1,0} & 1 & 0 & t_{1,3} & t_{1,4}\\
        t_{2,0} & 0 & 1 & t_{2,3} & t_{2,4}\\
    \end{bmatrix}.
\end{equation}
Here the rows are indexed by $(0,1,2)$ and columns are indexed by $(0,1,2,3,4)$.
By Remark \ref{re:paramflag} we see that this is a parametrization of the partial flag variety $\mathcal{F}l(1,2,3,4)$.
\end{example}

\begin{theorem}\label{thm:flag}
    When $\sigma = \{(1,1) ,(0,1), (1,0)\}$
    the truncation variety $V_\sigma$ is the flag variety $\mathcal{F}l(d - 1,d, d + 1,n)$.
\end{theorem}

\begin{proof}
    We recall that the truncated cluster operator $T(t_\sigma)$ is an element in the Fermi-Dirac algebra of the form
    $$
    T(t_\sigma) = T(t_{(1,0)}) +  T(t_{(1,1)}) + T(t_{(0,1)}) = \sum_{1 \le i \le d} t_{i,0} a_i + \sum_{1 \le i \le d < b \le n} t_{i,b} a_b^\dagger a_i + \sum_{d < b \le n} t_{0,b}a_b^\dagger.
    $$
    From Remark \ref{re:exteriorsofFD} we see that the operator $T(t_{(1,1)})$ commutes with the sum $T(t_{(1,0)}) + T(t_{(0,1)})$. Additionally, the sum is nilpotent of order $2$.
    Therefore, the truncation variety $V_{\sigma}$ is parameterized by the map $\mathcal{V}_\sigma \to \mathcal{F}$ where
    $$
    t_\sigma \mapsto \exp{T(t_{(1,1)})}e_{[d]} + T(t_{(0,1)})\exp{T(t_{(1,1)})}e_{[d]} + T(t_{(1,0)})\exp{T(t_{(1,1)})}e_{[d]}.
    $$
    The first summand $\exp{T(t_{(1,1)})}e_{[d]}$ is an element of the $d$th exterior power $\wedge^{d} \RR^n$ and by \cite[Theorem~3.5]{FSS} it parameterizes the Grassmannian $\operatorname{Gr}(d,n)$ in its Plücker embedding. We get that $\exp{T(t_{(1,1)})}e_{[d]} = t_1 \wedge \cdots \wedge t_d$ where 
    $
    t_i = (I + T(t_{(1,1)}))e_i
    $.
   The second summand is an element in the $d + 1$st exterior power $\wedge^{d + 1} \RR^n$ and 
    $$
    T(t_{(0,1)})\exp{T(t_{(1,1)})}e_{[d]} = T(t_{(0,1)}) t_1 \wedge \cdots \wedge t_d = t_{(0,1)} \wedge t_1 \wedge \cdots \wedge t_d.
    $$
    Here we set
    $t_{(0,1)} = T(t_{(0,1)}) 1
    $.
    Finally we note that the third summand is an element in the $(d - 1)$st exterior power $\wedge^{d - 1} \RR^n$. Dual to before we see that
    $$
    T(t_{(1,0)})\exp{T(t_{(1,1)})}e_{[d]} = T(t_{(1,0)}) t_1 \wedge \cdots \wedge t_d = t_{(1,0)} \lrcorner \, t_1 \wedge \cdots \wedge t_d.
    $$
    Here we set 
    $
    t_{(1,0)} = T(t_{(1,0)})^\dagger 1
    $.
    The coordinates of this element in $\wedge^{d - 1} \RR^n$ turn out to be $d \times d$ minors of the of a matrix
    $$
    T = \begin{bmatrix}
        0 & 0 & 0 & \cdots & 0 & t_{0,d + 1} & \cdots & t_{0,n}\\
        t_{1,0} & 1 & 0 & \cdots & 0 & t_{1,d + 1} & \cdots & t_{1,n}\\
        t_{2,0} & 0 & 1 & \cdots & 0 & t_{2,d + 1} & \cdots & t_{2,n}\\
        \vdots & \vdots & \vdots & \ddots & \vdots & \vdots & \ddots & \vdots\\
        t_{d,0} & 0 & 0 & \cdots & 1 & t_{d,d + 1} & \cdots & t_{d,n}\\
    \end{bmatrix}
    $$
    The rows of matrix $T$ are indexed by $(0,1, 2, \dots, d)$ and the columns are indexed by $(0, 1, 2, \dots, n)$.
    In fact, the coordinates of the whole parametrization are minors of $T$. A coordinate indexed by $J$, where $|J| = d + 1$, is equal to the maximal minor of $M$ taking columns $J$. A coordinate indexed by $J$ where $|J| = d$ is equal to the minor of $M$ taking the last $d$ rows and columns $J$. If $|J| = d - 1$ we do as before by adding the $0$th column. By Remark (\ref{re:paramflag}) this is a parametrization of $\mathcal{F}l(d - 1, d, d + 1,n)$.
\end{proof}

The theorem above is analogous to Thouless’ theorem and to \cite[Theorem~3.5]{FSS} in the context of ionization and electron attachment. We explain this connection for EOM-CC in the following remark.

\begin{remark}[Flag varieties in EOM-CC]
    We notice that when $\sigma = \{(1,1), (1,0) \}$ the truncation variety $V_{\sigma}$ is the flag variety $\mathcal{F}\ell (d - 1, d, n)$. From Theorem \ref{thm:eomcc} we see that the truncated EOM-CC ionized quantum states live in the Grassmannian
    $$
    \operatorname{Gr}(d - 1, n) = \mathcal{F}\ell (d - 1, d, n) \cap \wedge^{d - 1} \RR^n.
    $$
    Dually the truncated EOM-CC electron attachment quantum states are elements in the Grassmannian $\operatorname{Gr}(d + 1,n)$.
\end{remark}

In the title of this section, we also name spinor varieties; and for the rest of the section, we will study them and show how they appear in our constructions. Let $V = E \oplus F$ be a $2n$-dimensional vector space, where $E$ and $F$ are maximal isotropic subspaces with respect to some non-degenerate quadratic form $q$. Since $V$ is even dimensional, an isotropic subspace is \textit{maximal} if and only if it is $n$-dimensional.
The variety of maximal isotropic subspaces in $V$ splits into two isomorphic components denoted $S_-$ and $S_+$. Without loss of generality, we can assume that $E \in S_+$. Two maximal isotropic subspaces are in the same irreducible component if the dimension of their intersection has the same parity as $n$. The \textit{spinor variety} for $V$ is defined to be the irreducible component $S_+$. Its dimension is 
$\dim(S_+) = \binom{n}{2}$
and its degree is the number of shifted standard Young tableaux of shape $(n, n - 1, \dots, 1)$ \cite{hiller1982combinatorics}. 
See Manivel's paper \cite[Section~2]{manivel2009spinor} for a more detailed discussion on spinor varieties.

\begin{remark}[A parametrization of $S_+$]
    We choose a basis $\{e_1, \dots, e_n, f_1, \dots, f_n\}$ for $V$, such that the $e_i$'s form a basis for $E$ and the $f_i$'s form a basis for $F$. 
    We look at an $n \times 2n$ matrix $M = [\, I_n \, |\, U\, ]$ where $U$ is some $n \times n$ skew-symmetric matrix. The matrix $M$ is written over the basis described above. One can check that the rows of $M$ span an $n$-dimensional isotropic subspace $E_U$. We also get that 
    $$
    \dim(E \cap E_U) = n - \operatorname{rank}(U).
    $$
    Since $U$ is skew-symmetric, its rank is an even number, and hence the parity of the dimension of the intersection is the same as $n$.
    Hence $E_U \in S_+$.
    The subspace $E_U$ can be parametrized by the Pfaffians of $U$. 
    For a generic skew-symmetric matrix $U$ the subspace $E_U$ is a generic maximal isotropic subspace in $S_+$, so the map
    $$
    \PP^{\binom{n}{2} - 1} \to \PP^{2^{n - 1}}, \quad U \mapsto \text{ $2k \times 2k$ Pfaffians of $U$ for all $1 \le k \le n/2$}
    $$
    defines a birational parametrization of $S_+$.
\end{remark}

\begin{example}[A small spinor variety]
    We set $d = 2$, $n = 4$ and level set $\sigma = \{(1,1), (0,2), (2,0)\}$. The variety $V_\sigma$ is parametrized by the even Plücker~coordinates
    \begin{align*}
        &\psi_\varnothing = t_{12,0}, \,\, \psi_{12} = 1, \,\, \psi_{13} = t_{2,3}, \,\, \psi_{23} = -t_{1,3}, \,\, \psi_{14} = t_{2,4}, \\
        &\psi_{24} = -t_{1,4}, \,\, \psi_{34} = -t_{1,4}t_{2,3}+t_{1,3}t_{2,4}+t_{12,0}t_{0,34}, \,\, \psi_{1234} = t_{34}.
    \end{align*}
    These are the Pfaffians of a skew-symmetric matrix
    $$
    \begin{bmatrix}
        0 & -t_{12,0} & -t_{1,3} & -t_{1,4}\\
        t_{12,0} & 0 & t_{2,3} & t_{2,4}\\
        t_{1,3} & -t_{2,3} & 0 & -t_{0,34}\\
        t_{1,4} & -t_{2,4} & t_{0,34} & 0
    \end{bmatrix}.
    $$
    Coordinate $\psi_I$ is equal to the Pfaffian whose rows and columns come from the symmetric difference $I \oplus [d]$.
    This is a parametrization of the spinor variety for the vector space $\CC^8$. It is of dimension
    $$
    \dim(V_\sigma) = \binom{2}{2} + \binom{2}{2} + \binom{2}{1}\cdot \binom{2}{1} = 6
    $$
    and degree $\deg(V_\sigma) = 2$. It is defined by a single quadric, the homogeneous inverse master polynomial
    $$
    \overline{x}_{34}(\psi) = \psi_{23}\psi_{14}-\psi_{13}\psi_{24}+\psi_{12}\psi_{34} - \psi_\varnothing\psi_{1234}
    $$
    and the odd Plücker coordinates $\psi_1, \psi_2, \psi_3,\psi_4,\psi_{123}, \psi_{124}, \psi_{134}, \psi_{234}$.

\end{example}

\begin{theorem}\label{thm:spin}
    The truncation variety for $\sigma = \{(1,1), (0,2), (2,0)\}$ is the spinor variety of a $2n$-dimensional vector space $\CC^{2n}$.
\end{theorem}

\begin{proof}
Look at the $n$-dimensional vector space $F = \operatorname{span}\{a_1, \dots a_d, a_{d + 1}^\dagger, \dots, a_n^\dagger\}$. 
The basis vectors generate the subalgebra $\mathcal{S} \cong \wedge \CC^n$ of the Fermi-Dirac algebra, defined in Remark \ref{re:exteriorsofFD}. The truncated cluster operator is an element in the second exterior power of $F$:
$$
T(t_\sigma) = \sum_{1 \le i < j \le d}t_{ij, 0} a_ia_j + \sum_{1 \le i \le d < b \le n}t_{i, b} a_b^\dagger a_i + \sum_{d < b < c \le n}t_{0, bc} a_c^\dagger a_b^\dagger \in \wedge^2 F.
$$
The set $\{ (-1)^{i - 1}a_i : 1 \le i \le d \} \cup \{ (-1)^{d}a_b^\dagger : d < b \le n \}$ also forms a basis for $F$. We rewrite the operator $T(t_\sigma)$ in terms of this basis: 

\begin{footnotesize}
\begin{equation*}
    \sum_{1 \le i < j \le d} \hspace{-1em} (-1)^{i + j}t_{ij, 0}(-1)^{i  - 1 + j - 1}a_ia_j + \sum_{1 \le i \le d < b \le n} \hspace{-1em}(-1)^{i + d}t_{i, b} (-1)^{i - 1 + d}a_ia_b^\dagger - \sum_{d < b < c \le n} \hspace{-1em} t_{0, bc} (-1)^{2d}a_b^\dagger a_c^\dagger.
\end{equation*}
\end{footnotesize}

\noindent We set $T$ as the $n \times n$ skew symmetric matrix corresponding to $T(t_\sigma)$ with respect to the signed basis of $F$.
The exponential of $T(t_\sigma)$ is the operator
$$
\exp(T(t_\sigma)) = \sum_{k = 0}^d \frac{1}{k!}T(t_\sigma)^{\wedge k} = \sum_{\substack{I \subseteq [n] \\ |I| = 2k }} (-1)^{\ell d + \sum_{i \in [d] \cap I} (i - 1)}\operatorname{Pf}(T_I) a_{i_{1}} \cdots a_{i_{\ell }} a_{i_{\ell + 1}}^\dagger \cdots a_{i_{2k}}^\dagger 
$$
where $\operatorname{Pf}(T_I)$ is the Pfaffian of $T$ using columns and rows from $I$. For a compatible definition of the Pfaffian see \cite[Section~2.7.4.]{landsberg2011tensors}.
The exponential map is defined as
\begin{align*}
    \exp{T(t_\sigma)}e_{[d]} &= \sum_{\substack{I \subseteq [n] \\ I \text{ even}}} (-1)^{\ell d + \sum_{i \in [d] \cap I} (i - 1)}\operatorname{Pf}(T_I) a_{i_{1}} \cdots a_{i_{\ell}} a_{i_{\ell + 1}}^\dagger \cdots a_{i_{2k}}^\dagger a_1^\dagger \cdots a_d^\dagger 1\\
    % &= \sum_{\substack{I \subseteq [n] \\ I \text{ even}}} (-1)^{\sum_{i \in [d] \cap I} (i - 1)} \operatorname{Pf}(T_I) a_{i_{1}} \cdots a_{i_{\ell}} a_1^\dagger \cdots a_d^\dagger a_{i_{\ell + 1}}^\dagger \cdots a_{i_{2k}}^\dagger e_\varnothing\\
    &= \sum_{\substack{I \subseteq [n] \\ I \text{ even}}} \operatorname{Pf}(T_I) e_{[d] \oplus I}.
\end{align*}
The reason we introduced the signed basis for $F$ is to cancel out the signs appearing when multiplying the exponential operator $\exp{T}(t_\sigma)$ with the reference state $e_{[d]}$. We now see that the truncation variety is parametrized by the Pfaffians of a skew symmetric matrix $T$. This parameterizes the spinor variety for vector space $\CC^{2n}$.

% We want to show that the coordinates of the  exponential parametrization $\mathcal{V}_\sigma \to \mathcal{F}$ are the Pfaffians of the skew-symmetric matrix
% $$
% \begin{bmatrix}
%     0 & t_{12,0} & -t_{13,0} & \cdots & \pm t_{1,n}\\
%     -t_{12,0} & 0 & t_{23,0} & \cdots & \mp t_{2,n}\\
%     t_{13,0} & -t_{23,0} & 0 & \cdots & \pm t_{3,n}\\
%     \vdots & \vdots & \vdots & \ddots & \vdots\\
%     \mp t_{1,n} & \pm t_{2,n} & \mp t_{3,n} & \cdots & 0
% \end{bmatrix}.
% $$
% The coordinate $\psi_I$ is equal to the Pfaffian indexed by the columns the symmetric difference $I \oplus [d]$. 
% The signs arise from the exterior algebra, we define the sign of element in $(i,j)$-th position where $i < j$ as
% $
% (-1)^{\min({i, d}) + \min({j - 1, d})}
% $.
    
\end{proof}

We notice the truncation variety $V_{\{(2,0), (1,1), (0,2)\}} \cong S_+$ does not depend on the number of base electrons $d$. It only depends on the number of orbitals $n$.

\section{The coupled cluster equations}\label{se:CCeqs}

The central problem in electronic structure theory is to solve the electronic Schrödinger equation
$$
H\psi = \lambda\psi, \quad \psi \in \mathcal{F}
$$
for an electronic system of interest. Here $H$ is a $2^n \times 2^n$ symmetric matrix associated with the electronic system, called the \textit{Hamiltonian}. The solutions to the Schrödinger equation describe the stable states of our system; the eigenvalues $\lambda$ are the energies and the eigenvectors $\psi$ are the wave functions representing the probability amplitude of finding an electron at a particular location. The lowest eigenvalue and the corresponding eigenvector describe the \textit{ground state}, which is the most stable arrangement of electrons in the system. The higher energy solutions are called \textit{excited states} and they require energy input, such as from heat, electricity or light, to be observed. Hence, all solutions of the equation can be of interest.

For the sake of this paper, we will assume $H$ is a generic symmetric matrix of size $2^n \times 2^n$. For large $n$, solving the Schrödinger equation is infeasible. \textit{Coupled cluster theory} describes a hierarchy of approximations for the solutions of the Schrödinger equation, using the truncation varieties introduced in Section \ref{se:TruncVars}. We will now describe the construction of these approximation systems for a given truncation.

With that being said, let $\sigma$ be a proper subset of the truncation grid $\mathcal{G}$. We define $\psi_\sigma$ to be the projection of vector $\psi \in \mathcal{F}$ to coordinates with excitation level in $\sigma$. This is the projection of quantum states in $\mathcal{F}$ onto the subspace 
$$
\mathcal{F}_\sigma = \operatorname{span} \{ e_J : J \subseteq [n] \text{ and } (a(J), c(J)) \in \sigma \cup \{(0,0)\}\} \subseteq \mathcal{F}.
$$
We approximate the Schrödinger equation by restricting to quantum states on the truncation variety and by relaxing the constraints of the eigenvalue equations; only the projections of $H\psi$ and $\psi$ onto $\mathcal{F}_\sigma$ need to be linearly dependent. In symbols:
\begin{equation}\label{eq:CCcompact}
    (H\psi)_\sigma = \lambda\psi_\sigma, \quad \psi \in V_\sigma.
\end{equation}
These are called the \textit{unlinked coupled cluster (CC) equations} truncated at level set $\sigma$. Two variants of the coupled cluster equations are commonly used in the theory, the linked and unlinked CC equations.
See \cite[Section~4]{faulstich2024coupled} for the definition of the linked CC equations. By \cite[Theorem~5.11]{FSS} the two formulations agree in all cases that appear in the computational chemistry literature, including CCS, CCD, CCSD and CCSDT \cite{faulstich2024coupled}. The unlinked equations are more elegant from an algebraic point of view and its algebraic degree is lower, explaining our~choice.
The coupled cluster equations for systems with a fixed number of electrons correspond to choosing level sets $\sigma$ along the diagonal of the truncation grid. These cases are covered in detail in \cite[Section~5]{FSS}. In the same spirit as before, this section extends the coupled cluster framework to the full Fock space, thereby allowing both ionization and electron attachment. 

The dimension of the truncation variety $V_\sigma$ is equal to the number of constraints imposed by the equations in (\ref{eq:CCcompact}). Hence, the number of solutions to the CC equations is finite for a generic matrix $H$. We call this number the \textit{CC degree} of truncation variety $V_\sigma$, denoted ${\rm CCdeg}_{d,n}(\sigma)$. By \cite[Theorem~5.2]{FSS} we obtain an upper bound:
\begin{equation}\label{eq:upperbound}
    {\rm CCdeg}_{d,n}(\sigma) \le (\dim(V_\sigma) + 1)\deg(V_\sigma).
\end{equation}
By \cite[Corollary~5.3]{FSS} this bound is tight when $V_\sigma$ is linear. In that case the CC equations are 
an eigenvalue problem for the submatrix $H_{\sigma,\sigma}$ and 
${\rm CCdeg}_{d,n}(\sigma) = \dim(V_\sigma) + 1$. So if $H$ is a real symmetric matrix and the truncation variety is linear then all solutions to the CC equations are real.

Let $V$ be a variety parametrized by the rational map $\gamma : \PP^m \dashrightarrow
 \PP^s$.  We define the \textit{graph} $\Gamma_\gamma$ of $\gamma$ as the Zariski closure of the set $\{(x, \gamma(x)) : x \in \PP^m \}$. This is a variety in the product space $\PP^{m} \times \PP^s$ of same dimension as $V$. A graph $\Gamma_\gamma$ admits a natural $\ZZ^2$-grading with respect to its grading in $\PP^m$ and $\PP^s$. The $i$th bidegree of $\Gamma_\gamma$ is defined as
$$
\delta_i(\Gamma_\gamma) = \#(\Gamma_\gamma \cap (L \times L'))
$$
where $L \subseteq \PP^m$ and $L' \subseteq \PP^s$ are general linear subspaces of dimension $i$ and $s - i$, respectively. We can also define $i$th bidegree of a graph $\Gamma_\gamma$ to be the degree of the image of the restriction $\gamma(L) \subseteq V$ to a general linear space $L \subseteq \PP^m$ of dimension $i$.
We notice that $\delta_i = 0$ for $i < 0$ and $i > \dim(V)$. Additionally we see that $\delta_0 = 1$ and $\delta_s = \deg(V)$. 
The \textit{total degree} of graph $\Gamma_\gamma$ is defined to be the sum of its bidegrees,
$
\operatorname{tdeg}(\Gamma_\gamma) = \sum_{i = 0}^s \delta_i(\Gamma_\sigma).
$ Another way to define the total degree is to take its cone in affine space $\mathbb{A}^{m + 1} \times \mathbb{A}^{k + 1} \cong \mathbb{A}^{m + k + 2}$ and projectivize, obtaining a projective variety in $\PP^{m + k + 1}$. The degree of this variety is the total degree of the graph.
In \cite[Theorem~1]{borovik2025coupled}, we show that the CC degree of the Grassmannian is equal to the total degree of the graph of the restricted exponential map $\mathcal{V}_{\{1\}} \to \mathcal{F}$. We want to generalize this result.

A pair $(m,\ell)$ has a partition $\{(m_1,\ell_1), \cdots (m_k,\ell_k)\}$ over the level set $\sigma \subsetneq \mathcal{G}$ if 
$$
(m,\ell) = \sum_{i = 1}^k (m_i, \ell_i)
$$
and $(m_i, \ell_i) \in \sigma$ for all $1 \le i \le k$. We call $k$ the length of the partition. We say that this partition is even if at most one of the parts $(m_i,\ell_i)$ has an odd level $m_i + \ell_i$. If the element $(m,\ell)$ has an even level, then a set partition is even if and only if all of the parts have an even level. If the element  $(m,\ell)$ has an odd level, then a set partition is even if and only if exactly one part has an odd level. 

\begin{theorem}\label{thm:degofgraph}
    We take level set $\sigma \subseteq \mathcal{G}$ such that
    no element in $\sigma$ has an even partition over $\sigma$ of length $> 1$.
    The CC degree of the truncation variety $V_\sigma$ is the total degree of the graph of the restricted map
    \begin{equation}\label{eq:Vparam}
        \mathcal{V}_\sigma \to \mathcal{F}, \quad t \mapsto \psi = \exp{(T_\sigma(t))}e_{[d]}.
    \end{equation}
\end{theorem}

\begin{proof}
    We look at the projection $\psi_\sigma$ of $\psi$ onto $\mathcal{F}_\sigma$. 
    Take an indexing set $J$ with excitation level $(a(J), c(J))$ in $\sigma$. 
    By Theorem \ref{thm:psiparam} the exponential parametrization maps $\psi_J$ to a signed sum of monomials $t_\pi = t_{\pi_1 \cap [d], \pi_1 \backslash [d]} \cdots t_{\pi_k \cap [d], \pi_k \backslash [d]} $, where $\pi$ is an even set partition of $[d] \oplus J$.
    In the truncated parametrization (\ref{eq:Vparam}), we set $t_{I,B} = 0$ if $(|I|, |B|) \notin \sigma$. Hence, we only need to sum over monomials $t_\pi$ such that $(a(\pi_i), c(\pi_i)) \in \sigma$ for all $i$. We note that at most one of the levels $(a(\pi_i), c(\pi_i))$ is odd and $(a(J), c(J)) = \sum_{i} (a(\pi_i), c(\pi_i))$. 
    Assuming the hypothesis of the theorem the only even partition of $(m,\ell)$ over $\sigma$ is the trivial one $\{(m,\ell)\}$ of length $1$. Hence the monomials $t_\pi$ vanish for all even set partitions $\pi$ except the trivial partition $\{[d] \oplus J\}$. Thus $\psi_J = \pm t_{[d] \backslash J, J \backslash [d]}$.
    Now we use the same argument as in the proof of \cite[Theorem~1]{borovik2025coupled} to show that the CC degree is the total degree of the graph~of~(\ref{eq:Vparam}).
\end{proof}

In particular, the CC degree of the flag variety $\mathcal{F}l(d - 1,d, d + 1,n)$ is the total degree of the graph of the parametrization defined in Remark \ref{re:paramflag}. Also the CC degree of the spinor variety $S_+$ is the total degree of the graph of the map
$$
\PP^{\binom{n}{2} - 1} \to \PP^{2^{n - 1} - 1}, \quad \text{skew sym.\ } T \mapsto \operatorname{Pf}(T_I) \text{ for all even } \, I \subseteq [n],
$$
mapping skew symmetric matrix $T$ to its sub Pfaffians.
The degrees of the flag variety and the spinor variety are calculated using Schubert calculus. This begs the question of whether we can use Schubert calculus to find the total degrees of their graphs and therefore explain the fifth row of Tables \ref{ta:CCdegflag} and \ref{ta:CCdegSpin}? 

In \cite{borovik2025coupled, feigin2024birational,feigin2025birational} toric degenerations are used to find the total degrees of the graphs of birational parameterizations of the Grassmannian and the flag variety. Another approach for finding the CC degree of the flag variety and the spinor variety could be to look at toric degenerations. We leave these question open, but in the following section we offer a numerical study of CC degrees of flag varieties and spinor varieties.

\section{Numerical solutions}\label{sec:numsols}

The truncation variety $V_\sigma$ is parametrized by the restriction of the exponential map in (\ref{eq:expparam}) to the subspace $\mathcal{V}_\sigma$:
$$
\mathcal{V}_\sigma \to \mathcal{F}, \quad t_\sigma \mapsto \psi = \exp(T(t_\sigma))e_{[d]}.
$$
This map is birational and hence its fibers are zero dimensional and of degree $1$. Therefore we can rewrite the coupled cluster equations as:
\begin{equation}\label{eq:CCeqssimple}
    ((H - \lambda I )\exp(T(t_\sigma))e_{[d]})_\sigma = 0.
\end{equation}
This is a square polynomial equation of size $\dim(V_\sigma) + 1$ with polynomials in variables $(\lambda, t_\sigma)$ of degree $\le d + 1$. 
The above formulation works well for numerically solving the CC equations. In our computations we use \texttt{HomotopyContinuation.jl} \cite{breiding2018homotopycontinuation}, a package in Julia for numerically solving polynomial systems using homotopy methods. For a detailed introduction into numerical homotopy methods and how they can be used to both solve generic and special CC equations, we refer to a paper by the author and Fabian Faulstich  \cite{sverrisdottir2024exploring}. Throughout this section, we use Julia version 1.11.4 and \texttt{HomotopyContinuation.jl} version 2.9.2. Computations were done on the MPI-MiS computer server using four 18-Core Intel Xeon E7-8867 v4 at 2.4 GHz (3072 GB RAM). Symbolic calculations are done in the \texttt{Macaulay2} software system \cite{M2} on a MacBook Pro with 16 cores and an Apple M4 Max chip.

\begin{example}[CC degrees of flag varieties]
    We look at the truncation variety $V_{\{(1,0), (1,1), (0,1)\}} \cong \mathcal{F}\ell(d -1,d, d+1, n)$ when $d = 2,3$ and $n = 4,5,6,7,8$.
   
    \begin{table}[ht]
        \centering
        \caption{CC systems corresponding to the flag varieties.}
        \begin{tabular}{ccccccc}
            $(d,n)$ & (2,4) & (2,5) & (2,6) & (3,6) & (2,7) & (3,7)\\
            $\dim$ & 8 & 11 & 14 & 15 & 17 & 19\\
            ${\rm degree}$ & 12 & 110 & 1\,274 & 4\,550 & 17\,136 & 271\,320 \\
            ${\rm mingens}$ & $[2 ,\, 10]$ & $[7 ,\, 50]$ & $[23, \, 175]$ & $[14 , \, 281]$ & $[65, \, 490]$ & $[37, \, 1148]$\\
            ${\rm CCdeg}_{d,n}$ & 74 & 713 & 8\,499 & 30\,070 & 116\,602 & 1\,821\,528 \\
            \# real & 18 & 51 & 151 & 332 & 503 &  3\,262 \\
            solve (sec) & 1 & 6 & 174 & 712 & 14\,638 & 398\,022 \\
            certify (sec) & 0 & 0 & 2 & 4 & 24 & 329 \\
        \end{tabular}
        \label{ta:CCdegflag}
    \end{table}
    
    \noindent The dimension of $V_\sigma$ is $d(n - d) + n$. The degree and the minimal generators of $V_\sigma$ are calculated using \texttt{Macaulay2} \cite{M2}. The fourth row labeled ``mingens'' gives us the number of minimal generators of degree $1$ and $2$. The linear generators are the variables indexed by sets $J$ where $|J| < d - 1$ or $|J| > d + 1$. We note that the flag variety $ \mathcal{F}\ell(d -1,d, d+1, n)$ is cut out by quadrics, so we do not have any higher-degree generators. The CC degrees are calculated by solving (\ref{eq:CCeqssimple}) for a generic matrix, using the \textit{monodromy} function in \texttt{HomotopyContinuation.jl} \cite{breiding2018homotopycontinuation}. The row ``$\#$ real'' lists the number of real solutions to the CC equations for one generic real Hamiltonian. The number of real solutions might vary, and the number listed is only from one case.
\end{example}

Another truncation variety of interest is the spinor variety, $S_+$ of $\CC^{2n}$. From Theorem \ref{thm:spin} we see that the spinor variety, as a truncation variety, only depends of the number of orbitals $n$ but not on the number of electrons $d$. 

\begin{example}[CC degrees of spinor varieties]
    We look at the truncation variety $V_{\{(2,0), (1,1), (0,2)\}} \cong S_+$, isomorphic to the spinor variety of $\CC^{2n}$, when $n = 4,5,6,7$.
    \begin{table}[ht]
        \centering
         \caption{CC systems corresponding to the spinor variety.}
        \begin{tabular}{cccccc}
            $n$ & 4 & 5 & 6 & 7\\
            $\dim$ & 6 & 10 & 15 & 21\\
            ${\rm degree}$ & 2 & 12 & 286 & 33\,592\\
            ${\rm mingens}$ & $[8, \, 1]$ & $[16, \,10]$ & $[32, \, 66]$ & $[64, \, 364]$\\
            $\text{CCdeg}_{d,n}$ & 13 & 98 & 2\,572 & 318\,118 \\
            \# real & 5 & 10 & 70 & 904 \\
            solve (sec) & 0 & 1 & 32 & 60\,504 \\
            certify (sec) & 0 & 0 & 1 & 55\\
        \end{tabular}
        \label{ta:CCdegSpin}
    \end{table}
    
    \noindent The dimension of $V_\sigma$ is $\binom{n}{2}$.  The linear generators of the truncation variety are the variables indexed by the odd set $J \subseteq [n]$. We note that the spinor variety is cut out by quadratics so we only have minimal generators of degree $1$ and $2$.
\end{example}

Similar tables can be made for the truncation varieties corresponding to CCD and CCSD. The calculations become large quickly as can be seen in the final example.

\begin{example}
    We look at FSCCSD for ionized electronic systems with $d$ base electrons in $n$ orbitals. Our level set is of the form $\sigma = \{(1,0), (2,1), (1,1), (2,2)\}$. The corresponding truncation variety is of dimension
    $$
    \dim(V_\sigma) = d + d(n - d) + \binom{d}{2}(n - d) + \binom{d}{2}\binom{n - d}{2}.
    $$
    By Theorem \ref{thm:linear} the truncation variety $V_\sigma$ is linear when $d = 2$.
    We now focus on the case when $d = 3$. The dimension of $V_\sigma$ is equal to 
    $$
    30, 45, 63, 84, 108, 135, 165, ... \quad \text{ for } n = 6,7,8,9,10,11,12,... \,.
    $$
    Solving square polynomial systems of these sizes is hard, even numerically. Using \texttt{HomotopyContinuation.jl} we are able to compute the CC degree of $V_\sigma$ when $n = 6$. In this case $V_\sigma$ has degree $43$ and thus (\ref{eq:upperbound}) provides us with an upper bound of $31 \cdot 43 = 1333$. Computations using \texttt{monodromy} show that the actual CC degree is 
    $$
    \operatorname{CCdeg}_{3,6}(\sigma) = 1195 < 31 \cdot 43 = 1333.
    $$
\end{example}

\bibliographystyle{siamplain}
\bibliography{bibtex}

\end{document}